\documentclass[12pt]{amsart}
\usepackage{amssymb,latexsym,amsmath, amscd}

\theoremstyle{plain}
\newtheorem{theorem}{Theorem}[subsection]
\newtheorem{lemma}[theorem]{Lemma}
\newtheorem{proposition}[theorem]{Proposition}
\newtheorem{corollary}[theorem]{Corollary}
\newtheorem{definition}[theorem]{Definition}

\newtheorem{ttheorem}{Theorem}[section]
\newtheorem{llemma}[ttheorem]{Lemma}
\newtheorem{pproposition}[ttheorem]{Proposition}

\theoremstyle{remark}

\newtheorem{notation}[theorem]{\textit{Notation}}
\newtheorem{remark}[theorem]{\textit{Remark}}

\newtheorem{nnotation}[ttheorem]{\textit{Notation}}
\newtheorem{rremark}[ttheorem]{\textit{Remark}}

\numberwithin{equation}{subsection}

\newcommand{\be}{\mathbf{e}}
\newcommand{\br}{\mathbf{r}}
\newcommand{\bs}{\mathbf{s}}
\newcommand{\barbs}{\overline{\bs}}
\newcommand{\bP}{\mathbf{P}}
\newcommand{\bx}{\mathbf{x}}
\newcommand{\by}{\mathbf{y}}
\newcommand{\hx}{\hat\bx}
\newcommand{\hy}{\hat\by}
\newcommand{\bbR}{\mathbb{R}}
\newcommand{\bbRn}{\bbR^n}
\newcommand{\bbRnp}{(\bbR_+)^n}
\newcommand{\bbC}{\mathbb{C}}
\newcommand{\bbCm}{\bbC^{m+1}}
\newcommand{\id}{\mathbb{I}}
\newcommand{\bbE}{\mathbb{E}}
\newcommand{\bbCPm}{\bbC\mathbb{P}^m}

\newcommand{\calH}{\mathcal{H}}
\newcommand{\calHm}{\calH_{m+1}}
\newcommand{\calHmz}{\calHm^0}
\newcommand{\orb}{\mathcal{O}}
\newcommand{\calSm}{\mathcal{S}_{m+1}}
\newcommand{\cm}{\mathcal{M}}

\newcommand{\mur}{\mu_\br}
\newcommand{\vL}{{\varLambda}}
\newcommand{\laij}{\lambda_{ij}}
\newcommand{\laik}{\lambda_{ik}}
\newcommand{\lakl}{\lambda_{k\ell}}
\newcommand{\lail}{\lambda_{il}}

\newcommand{\phiij}{\phi_{ij}}

\newcommand{\thij}{\theta_{ij}}

\newcommand{\Mr}{M_\br}
\newcommand{\Mtw}{\widetilde M}
\newcommand{\Mrtw}{\Mtw_\br}

\newcommand{\Nrtw}{\widetilde N_\br}
\newcommand{\sigr}{\Sigma_\br}
\newcommand{\sigrtw}{\widetilde\Sigma_\br}

\newcommand{\wtens}[1]{w_{#1}\otimes w_{#1}^*}
\newcommand{\ei}{\sqrt{-1}\,}
\newcommand{\pair}[2]{\langle{#1},{#2}\rangle}

\newcommand{\Tr}{\mathop{\rm Tr}\nolimits}
\newcommand{\Ad}{\mathop{\rm Ad}\nolimits}
\newcommand{\Real}{\mathop{\rm Re}\nolimits}
\newcommand{\Imag}{\mathop{\rm Im}\nolimits}
\newcommand{\spann}{\mathop{\rm span}\nolimits}
\newcommand{\diag}{\mathop{\rm diag}\nolimits}
\newcommand{\Pol}{\mathop{\rm Pol}\nolimits}
\newcommand{\CPol}{\text{CPol}}
\newcommand{\Imtr}{\mathop{\rm Im\ Tr}\nolimits}

\newcommand{\SUm}{\,\text{SU(m+1)}}
\newcommand{\sutwo}{\mathfrak{su}(2)}
\newcommand{\um}{\mathfrak{u}(m+1)}
\newcommand{\essum}{\mathfrak{su}(m+1)}
\newcommand{\Um}{{\text{U(m+1)}}}
\newcommand{\Un}{{\text{U(n)}}}

\newcommand{\GLm}{{\text{GL(m+1,$\bbC$)}}}
\newcommand{\GLn}{{\text{GL(n,$\bbC$)}}}

\begin{document}

\title[Bending flows]
{Bending flows for sums of rank one matrices}
\author{Hermann Flaschka and John Millson}
\address{Department of Mathematics\\
         The University of Arizona\\
         Tucson, AZ 85721}
         \email{flaschka@math.arizona.edu}

\address{Department of Mathematics\\
         University of Maryland\\
         College Park, MD 20742}
         \email{jjm@math.umd.edu}


\maketitle

\tableofcontents

\section{Introduction}

The aim of this paper is to generalize results of
Kapovich and Millson \cite{KapovichMillson2} and Klyachko \cite{Klyachko1} on
the moduli
space $\Mr$ of (closed) polygons in $\bbR^3$ with prescribed sidelengths
$\br=(r_1,\ldots,r_n)$, polygons related by a Euclidean motion
being identified. They showed that $\Mr$ is a (possibly singular)
symplectic manifold, and introduced a class of commuting Hamiltonian
flows, the so-called  \emph{bending flows}. These flows
bend the polygon about the diagonals emanating from one fixed vertex.
The part of the polygon to one side of the diagonal does not move,
while the other part rotates at constant speed. The lengths of
the diagonals are action variables
which generate the bending flows; the conjugate angle
variables are the dihedral angles between the fixed and the moving
parts.

We generalize this picture by replacing vectors in $\bbR^3$ by
positive semidefinite rank-one Hermitean matrices. These have
the form $e=r\wtens{}$, where $r>0$ and $ w$ is a unit
vector in $\bbCm$. Explicitly,
$e:v\mapsto r(v,w)w$ where $(\ ,\ )$ is the standard positive definite
Hermitean form on $\bbCm$. The edges of a polygon will be $e_i=r_i\wtens{i},
i=1,\ldots, n$, $r_i>0$ fixed,
and \emph{closed} will mean \emph{closed up
to a multiple of the identity},
\[
e_1+\cdots+e_n=\vL\id.
\]
Equality of traces forces $\vL=(r_1+\cdots+r_n)/(m+1)$. Polygons
are identified if they are related by simultaneous rotation of
the sides by an element of $\Um$. Since $\wtens{}$ is unchanged if
$w$ is multiplied by $\exp(\ei\theta)$, we may think of
an edge $r\wtens{}$ as a weighted point in the projective space
$\bbCPm$, and of a polygon as a weighted configuration of points
in $\bbCPm$.

The paper has three parts.
\begin{enumerate}
\item\  A study of the moduli space $\Mr$. Criteria for nonemptiness
and nonsingularity of $\Mr$ are found.
\item\ A generalization of bending flows and their action-angle
coordinates.
\item\  A relation between the bending flows and representations
of $\Um$; this is reminiscent of geometric quantization.
\end{enumerate}
This is the logical progression of the material, but (1) stands
alone, and (2), about the bending flows, can be read independently
of the rest. We give a brief outline of each part.

\vspace{.07in}

If follows from general results of \cite{Klyachko2} 
(however we sketch a direct proof in what follows) that
$\Mr$ is nonempty if, and only if, the side lengths $\br$
satisfy the \emph{generalized triangle inequalities},
\[
mr_i\le r_1+\cdots+\hat r_i+\cdots r_n,\quad 1\le i \le n.
\]

This defines a convex cone $C(n,m+1)$ in $\bbRn$.
When $n=3$, the generalized triangle
inequalities are the usual conditions for $r_1,r_2,r_3$
to be the sides of a Euclidean triangle, namely
\[
r_1\le r_2+r_3,\quad r_2\le r_1+r_3,\quad r_3\le r_1+r_2.
\]

The proof of necessity is elementary, and is given in \S 3.1.
Sufficiency is deeper. Only an outline of the argument is given
in \S\S\ 3.2 and 3.3; for details the reader is referred to the
literature.
It follows from standard results, see for example  \cite{MumfordFogartyKirwan}, 
Theorem 8.3, that
$\Mr$ is canonically homeomorphic
to a weighted (by $\br$) complex analytic quotient of
the $n$-fold product $\Pi_1^n\, \bbCPm$.  For $m=1$, i.\ e.\ for
spatial polygons, the connection between weighted analytic
quotients of $\mathbb{CP}^1$ and $\Mr$ was found independently
in \cite{KapovichMillson2} and \cite{Klyachko2}.
The weighted quotient is nonempty if and only if there exists
a weighted semistable configuration of points on $\bbCPm$.
If $\br$ satisfies the strong triangle inequalities, then any $n$-tuple
of points on $\mathbb{CP}^m$ in general position (these always exist) is
semistable for
the weights $\br$.

In the Euclidean case, it was shown in \cite{KapovichMillson2} that
$\Mr$ is smooth if there is no polygon contained in
a line. The corresponding statement in our case is that $\Mr$ is
smooth if there is no \emph{decomposable} polygon.
A polygon is decomposable, roughly speaking, if there are closed
subpolygons contained in orthogonal subspaces of $\bbCm$. The
side lengths $\br$ for which a decomposable polygon exists lie on
certain hyperplane sections of the cone $C(n,m+1)$ called \emph{walls}.
A connected component of the complement of the walls is a \emph{chamber}.
We prove that $\Mr$ and $M_{\br'}$ are diffeomorphic if $\br,\br'$
lie in the same chamber. The topology of $\Mr$ will change
as $\br$ crosses a wall(see for example \cite{Hu} and \cite{Goldin}).

\vspace{.07in}

Turning now to the bending flows, we remark first that the action
variables, or Hamiltonians generating the ``bending'', are a
natural generalization of the Euclidean case. There, one has
$\bbR^3\equiv \sutwo$; a diagonal $A$ of a polygon
is identified with
an element $\hat A\in \sutwo$, and the length $\|A\|$ is just the
positive eigenvalue of $\hat A$. In our generalization, the
bending Hamiltonians are also the eigenvalues, $\laij,j=1,\ldots,m+1$,
of the diagonals
$A_i=e_1+\cdots+e_{i+1}$, and they generate $2\pi$-periodic flows.
These flows again leave part of the polygon fixed, and conjugate
the other part by $\exp(\ei tE_{ij})$, where $E_{ij}$ is
the spectral projection for $\laij$. This is a kind of bending
with ``internal degrees of freedom''.
Enough of the $\laij$ are functionally independent to give
action coordinates (only on a dense open set, however). The
eigenvalues of $A_{i+1}$ and $A_i$ interlace; this is a
direct consequence of the Weinstein-Aronszajn formula from
perturbation theory. The interlacing property can be pictured
by a triangular organization of the $\laij$, starting with
the eigenvalue $r_1$ of $e_1=r_1\wtens{1}$ and working up:
\begin{equation}
\begin{matrix}
   \lambda_{21}&            &\lambda_{22}&            &\lambda_{23}\\
               &\lambda_{11}&            &\lambda_{12}&            \\
               &            &     r_1    &            &
\tag{*}
\end{matrix}
\end{equation}
This is called a \emph{Gel'fand-Tsetlin} pattern.
As long as these inequalities are respected, the $\laij$ can be
prescribed arbitrarily.
Gel'fand-Tsetlin patterns were introduced in the
Euclidean context by Hausmann and Knutson \cite{HausmannKnutson}.
They showed,
and this extends easily to our setting, that the bending flows
and the Gel'fand-Tsetlin flows of Guillemin and Sternberg
\cite{GuilleminSternberg} are dual to each other (via the
Gel'fand-MacPherson duality \cite{Gel'fandMacPherson}).

The angle variables also extend the Euclidean case (dihedral angles),
but in a rather more subtle way. For vectors
$w,x,y,z$ in $\bbCm$, one defines the
\emph{four-point function} \cite{BerceanuSchlichenmaier}
\[
F(w,x,y,z)=(w,x)(x,y)(y,z)(z,w).
\]
It is independent of the phases of its arguments, and
so $\arg F$ is well defined on $\bbCPm$. The arguments of $F$ in our
setting will be the $w_i$ that define the edges $e_i=r\wtens{i}$, and
eigenvectors $u_{ij}$ corresponding to $\laij$. These are only
defined up to phase, but using
$\arg F(w_{i+1},u_{ij},w_{i+2},u_{i,j+1})$ we get
global angle variables (on a dense open set again). In the
Euclidean case, this amounts to a rather complicated way
of expressing a dihedral angle. The proof of the
Poisson bracket relations, $\{$angle , angle$\}=0$, etc.,
takes up all of \S 7.

\vspace{.07in}

The connection between bending flows and representation theory
follows the ideas of Guillemin and Sternberg \cite{GuilleminSternberg}.
In Bohr-Sommerfeld quantization, one asks that the action
variables take on integral values. If the $r_i$, $\laij$, and
$(\sum_i r_i)/(m+1)$ are integers, the interlacing property
of the $\laij$ reproduces the Pieri formula for the decomposition
of tensor products of symmetric powers
\begin{equation}
\bigotimes \mathcal{S}^{r_i}(\bbCm) \tag{**}
\end{equation}
of the basic representation of $\Um$ on $\bbCm$.

The Hausmann-Knutson duality also has a representation-theoretic meaning.
Gel'fand-Tsetlin patterns were invented
to parameterize bases for vector spaces carrying representations
of unitary (and general linear) groups. The possible patterns
(*) built from (integer)
eigenvalues of successive
diagonals index vectors of weight $\br=(r_1,\ldots,r_n)$ in the
(Grassmannian) representation of $\Un$ with highest weight
\begin{equation}
(\underbrace{\vL,\ldots,\vL}_{m+1},0,\dots,0). \tag{***}
\end{equation}
We conclude that
the multiplicity of the one-dimensional representation
$\det{}^\vL$ in the tensor product (**) of $\Um$ representations
equals the
multiplicity of the weight $\br$ in the representation (***) of $\Un$.

It must be said that these results are based on counting lattice
points in convex polytopes and comparing their number with
multiplicities known in representation theory. We do not construct
actual representation spaces by any quantization method.

\vspace{.07in}

It is our hope that there are analogous results for all symplectic
quotients
of products of flag manifolds.
In general for such products, one can find integrable systems
that reduce to ours in the
case of projective space, but it appears very hard to construct
an explicit
family of Hamiltonians {\em with periodic flows}, i.e. action variables.
If such a construction could be carried out and the associated momentum
polytope could be
computed, then by counting lattice
points in the momentum polytope one could find information on decomposing
tensor products of irreducible representations. Many deep connections
are now known between tensor product decompositions and convex
polyhedra; these, however, do not seem to arise as
images of momentum mappings. One of the main motivations
for our paper is that the special case treated here is
probably the only case
where everything can be worked out with simple explicit formulas.

We conclude by noting that the spaces $M_{\br}$ studied in this paper  were the
subject of the book, \cite{DolgachevOrtland}. The study of the spaces
$M_{\br}$ in \cite{DolgachevOrtland} was from the point of view of
algebraic geometry and combinatorics, necessitating the restriction to the
case in which the $r_i$'s were integral (moreover the authors assumed
that all the $r_i$'s were equal).  There appear to be interesting relations
between our work and theirs.

\vspace{.15in}

\textbf{Acknowledgements}. We thank the referee for detailed and
helpful comments, and especially for pointing out an error
in the original version of this paper. In the appendix we answer
a question posed by the referee concerning the relation between the duality of
integrable systems of this paper and those of \cite{AdamsHarnadHurtubise}.
We would also like to thank Ron Donagi for some helpful conversations
about Hitchin Hamiltonians. He suggested to us that if the bending Hamiltonians
belonged to the Hitchin system then the underlying curve would have to 
degenerate (see the appendix).  Finally, we thank Ben Howard for drawing our attention to
confusion caused by a missing factor of $2$.
The second author was supported in part by NSF grant DMS  01-04006.

\section{The moduli space of polygons in $\calHm$}

In this section, we collect the notation used throughout, and in particular,
introduce the moduli space of polygons with which we will be concerned.

\subsection{Coadjoint orbits}

\begin{enumerate}
\item Let $\calHm$ be the vector space of $(m+1)\times(m+1)$ Hermitean matrices. We
identify it with the dual of
the Lie algebra $\um$ via the pairing
$\pair{\xi}{X}=\Imtr \xi X$, for $\xi\in\um, X\in\calHm$.

\item $\calHmz=\{X\in\calHm\mid\Tr X=0\}$. It is the dual of $\essum$.

\item The gradient $\nabla f(X)\in\um$ of a smooth
function $f$ on $\calHm$ is defined by
\[
 \Imtr(\nabla f(X)Y) = \frac{d}{dt}\bigg|_{t=0} f(X + t Y),
 \quad\mbox{for\ all\ }Y \in \calHm.
\]

\item  A $\Um$-orbit $\orb \subset \calHm$  carries the
Kostant-Kirillov symplectic form $\omega_{KK}$ defined by
\[
\omega_{KK}(X)([\xi,X],[\eta,X])=\Imtr(X[\xi,\eta]).
\]
The Lie-Poisson bracket is
\begin{equation}\label{bracket}
\{f,g\}(X)=\Imtr(X\,[\nabla f(X),\nabla g(X)]),
\end{equation}
and Hamilton's equations have the form
\begin{equation}\label{Ham-equn}
\dot X=[\nabla f(X),X].
\end{equation}
In \S7.1, we use the identification between $\calH_2^0$ with bracket
\eqref{bracket}, and Euclidean space $\bbR^3$ with its standard Poisson
bracket.

\item Let $w\in\bbCm$ be a unit vector. Define $\wtens{}\in\calHm$
by $\wtens{}(v)=(v,w)w$; it is a rank one projection. The matrices
$r\wtens{}, r>0$, form an orbit $\orb_r$ of $\Um$. They will
be ``edges" of polygons, and are denoted by the letter $e$.
Given $e\in\orb_r$, the unit vector $w$ is determined up to
multiplication by a complex number of modulus one.  Hence
$\orb_r$ is diffeomorphic to $\bbCPm$. As symplectic manifolds, they
are related by $\omega_{KK}=2r\omega_{FS}$, where $\omega_{FS}$ is the
Fubini-Study form on $\bbCPm$.
\end{enumerate}

\begin{remark} For completeness, we verify the last assertion about
the symplectic forms. The Fubini-Study metric is
the $U(m+1)$-invariant K\"ahler metric normalized
so that the holomorphic sectional curvature is $4$.
To determine the scalar multiple relating two $U(m+1)$
invariant symplectic forms on $\bbCPm$ it suffices to compute the period of each
over a (complex linearly) embedded projective line $\mathbb{CP}^1 \subset \mathbb{CP}^m$.
For the Fubini-Study metric, a projective line has curvature $4$, and is
therefore isometric to a sphere of radius $1/2$ and area $\pi$.
We may obtain a such a $\mathbb{CP}^1$
by embedding $\mathcal{H}_2^0 \subset \calHm$ into the principal $2\times 2$
block.  Thanks to the identification of $\mathcal{H}_2^0$
with $\bbR^3$, it now suffices to relate the Kostant-Kirillov form on
an orbit of $\sutwo$ on $\bbR^3$, i.\ e.\ on a sphere
$S^2_r$ of radius r, to the usual area
form $dA$. The calculation in \cite[p. 412]{MR} shows that $\omega_{KK}(S^2_r)=
dA/2r$. The symplectic area is $2\pi r$, which agrees with the area
given by the scaled Fubini-Study form,
$2r\omega_{FS}$, on the embedded $\mathbb{CP}^1\subset\orb_r$.
\end{remark}

\subsection{The space of closed polygons}

Let $\br =(r_1,r_2,\ldots,r_n)$ be an $n$-tuple of positive numbers.
We define
a {\em (closed) polygon with side-lengths $\br$} to be an $n$-tuple
$\be  = (e_1,e_2,\ldots,e_n)$ such that for all $i,1 \le i \le n$ we have
\begin{itemize}
\item{a)} $e_i\in \orb_{r_i}$,
\item{b)} $\sum_1^n e_i = \vL \id$.
\end{itemize}
Note that
$\vL= \frac{1}{m+1}\sum_1^{m+1} r_i$ follows from equality of traces in b).
We call the matrices $e_i$ the {\em edges} of the polygon
$\be$ and
$r_i$ the {\em length} of the edge $e_i$.
Condition (b) says that the polygon $\be$ is closed, modulo
the center of $\calHm$.

\begin{enumerate}
\item When $\br$ is given, $\vL$ always stands for
$\frac{1}{m+1}\sum r_i$. Sometimes the notation $\vL_\br$ is
used to emphasize the dependence of $\vL$ on $\br$.

\item Given $\br$, define $\Nrtw$ to be the product symplectic manifold
$\Pi_1^n\, \orb _{r_i}$.  The diagonal action of $\Um$ on $\Nrtw$ is
Hamiltonian with momentum map $\mur$ given by
$$\mur(\be)=\sum_1^n e_i.$$
We refer to elements of $\Nrtw$ as \emph{linkages}; they may or may
not be closed.

\item Given $\br$, let
$$\Mrtw=\mur^{-1}(\vL\id)=\{\be\in\Nrtw \mid \sum_{i=1}^n e_i=\vL\id\}.$$
The elements of $\Mrtw$ are the \emph{polygons}, or \emph{closed
polygons} for emphasis. The unitary group acts diagonally on
$\Mrtw$.  We let ${\bf i}:\Mrtw \to \Nrtw$ be the inclusion.

\item Finally, we define the moduli space, $\Mr$, of polygons
(with side-lengths $\br$) to be the quotient of $\Mrtw$ by the diagonal
action of $\Um$.

\end{enumerate}

Because the stabilizer of the scalar matrix $\vL\id$ is all of $\Um$,
we obtain

\begin{lemma}
$\Mr$ is the symplectic quotient of $\Nrtw$ corresponding to the (one-point)
orbit $\vL \id \in \calHm$ under the diagonal action of $\Um$.
\end{lemma}

\section{Nonemptiness of the moduli spaces}

A simple set of inequalities on the side-lengths $r_i$ is necessary
and sufficient for the moduli space to be nonempty. The elementary
proof of necessity is given first. Sufficiency is deeper,
and is based on the interpretation of polygons as weighted sets of
points in projective space. For sake of completeness, that argument
will be summarized in the Subsection~3.2, with references to literature
where details may be found.

\subsection{Necessity of the Triangle Inequalities}

\begin{theorem}\label{nonempty}
The moduli space $\Mr$ is nonempty if and only if $\br$ satisfies the system
of inequalities
$$mr_i \le r_1 + r_2 + \cdots + \widehat{r_i} + \cdots +r_n, \
1 \leq i \leq n.$$
\end{theorem}
Here $\widehat{r_i}$ means that $r_i$ has been omitted in the summation.

\begin{remark}
We will call this system of inequalities (together with the inequalities
$r_i \geq 0, \ 1 \leq i \leq n)$ the {\em strong triangle inequalities of
weight} $m$. When $m=1$, they give the familiar conditions $r_1\le
r_2+r_3$, $r_2\le r_1+r_3$, $r_3\le r_1+r_2$ on the side lengths of
a planar triangle.
We will omit reference to the weight $m$ when it is clear from the context.
Note that if we define $\rho = \sum_i r_i$ then the $i$-th inequality
above is equivalent to
\begin{equation}\label{rhotriangle}
r_i \leq \frac{1}{m+1}\, \rho.
\end{equation}
\end{remark}

We now turn to the proof of the necessity of the triangle inequalities.

\begin{definition}
Let $X\in \calHm^0$. Say that $X$ is {\em maximally singular} if $X$
is conjugate to a diagonal matrix with eigenvalues
$(r,-\frac{r}{m},\ldots,-\frac{r}{m})$.
We note that the orbit $\orb_r^0$ under $\Um$ of such an $X$ is
the projection onto tracefree matrices of the orbit $\orb_r$
through $\diag(r,0,\ldots,0)$.
\end{definition}

\begin{lemma}\label{maxsing}
Suppose $X_1,X_2 \in \calHm ^0$ are distinct, maximally singular, and satisfy
$\Tr(X_j^2)=1$. Then
$\Tr (X_1X_2) \geq-1/m$, with equality if and only if $X_1$ and $X_2$ commute.
\end{lemma}

\begin{proof}
We may write
$$X_j=\sqrt{\frac{m+1}{m}}(\wtens{j}-\frac{1}{m+1}\id),$$
where $\|w_j\|=1,\,j=1,2$.

Then
\begin{align*}
\Tr X_1X_2 = &\frac{m+1}{m} \Tr[(\wtens{1} - \frac{1}{m+1} \id)
(\wtens{2})]\\
  =& \frac{m+1}{m}[|(w_1,w_2)|^2 - \frac{1}{m+1}] \geq
 \frac{m+1}{m}
\cdot - \frac{1}{m+1}\\
 =&-\frac{1}{m}.
\end{align*}
Clearly we have equality if and only if $(w_1,w_2) = 0$ if and only if
$X_1$ and $X_2$ commute.
\end{proof}

\begin{proposition}
Suppose that $\Mr$ is nonempty. Then $\br$ satisfies the strong triangle
inequalities of weight $m$.
\end{proposition}

\begin{proof}
Choose $\be \in \Mrtw$. Then $e_1+\cdots+e_n=\vL\id$ is
equivalent to $r_1X_1+\cdots+r_nX_n=0$, where the matrices
$$X_j=\sqrt{\frac{m+1}{m}}
(\wtens{j}-\frac{1}{m+1}\id)$$ satisfy the hypotheses of Lemma
\ref{maxsing}.
Alternatively,
$$r_iX_i=-r_1X_1-\cdots-\widehat{r_iX_i}-\cdots-r_nX_n.$$
Multiply each side by $X_i$ and take the trace to obtain
$$r_i^2=-\sum_{j(\ne i)} r_ir_j\Tr(X_j X_i) \leq\frac{1}{m}\sum_{j(\ne i)}
r_i r_j.$$
Now divide both sides by $r_i$ to obtain the result.
\end{proof}

The generalized triangle inequalities define a cone in $(\bbR_+)^n$:
\begin{definition}\label{hypersimplex}
$$C(n,m+1)=\{\br\in(\bbR_+)^n \mid \Mr\ne\emptyset\}.$$
\end{definition}
The intersection of $C(n,m+1)$
with the hyperplane $\sum r_i=m+1$ is known in the literature as
the {\em hypersimplex}. Side lengths $\br$ for which the moduli space
$M_\br$ is singular will be shown to lie on certain hyperplane sections
of $C(n,m+1)$.

We next give the outline of the proof that $M_\br$ is nonempty when
$\br\in C(n,m+1)$. The method is based on the identification of
an edge $r_i \wtens{i}$ with a point
$\bbC w_i\in\bbCPm$ of weight $r_i$. The nonemptiness of the
moduli space of such weighted points follows from a comprehensive
general theory.
\subsection{Semistability and Sufficiency of the Triangle Inequalities}

\subsubsection{Analytic quotients and symplectic quotients}

In \cite{Sjamaar} and \cite{HeinznerLoose}, the authors constructed
the analytic
quotient of a (not necessarily projective) compact K\"ahler manifold
$M$ by the action
of a complex reductive group $G$. It is assumed that some maximal compact
subgroup $K\subset G$ acts in a Hamiltonian fashion on $M$ with momentum
map $\mu$. In their theory, a point $m\in M$ is defined to be
{\em semistable} if the closure of the orbit $G\cdot m$ intersects the
subset
$\mu^{-1}(0)$ of $M$. The set of semistable points is denoted by $M^{sst}$;
it is open in $M$. A point of $M$ is defined to be {\em nice semistable} if
the orbit itself intersects $\mu^{-1}(0)$.
Define an equivalence relation, called {\em extended orbit equivalence},
by declaring two points to be related if their orbit closures intersect.
(That this is indeed an equivalence relation follows from a theorem
asserting that each equivalence class of semistable points contains a
unique nice semistable orbit). We emphasize that the notion of semistability
{\em depends on the symplectic structure}, i.e. on $\br$ in our case.

The {\em analytic quotient} of $M$ by $G$, denoted $M//G$,  is then defined
to be the quotient of $M^{sst}$
by extended orbit equivalence.  It is a hard theorem of \cite{Sjamaar} and
\cite{HeinznerLoose} that the resulting quotient topological space is
Hausdorff and
compact, and in fact has a canonical structure of a complex analytic space.

Since, by definition,
any point in $\mu^{-1}(0)$ is (nice) semistable, there is an induced map from
the symplectic quotient $\mu^{-1}(0)/K$  to the analytic quotient.
\cite{Sjamaar} and \cite{HeinznerLoose} prove
that this map is a homeomorphism.

These results were proved earlier for
smooth quotients and when
the quotient has only orbifold singularities (i.e. the stabilizer
of every $x \in \mu^{-1}(0)$ is finite) in \cite{Kirwan}, Theorem 7.5,
and also for general point stabilizers if $M$ is a smooth complex
projective variety
and the symplectic
form represents the (dual of the) hyperplane section class in
\cite{Kirwan},
Remark
8.14. In our setting, this amounts to assuming that either
$\br$ is not on a wall or that $\br$ is integral (it is probable that
some weakening of this condition will still result in an
integral K\"ahler class). In this case the
analytic quotient $M//G$ is a complex projective
variety.

We have seen that our space $M_{\br}$ is a symplectic quotient of
${(\bbCPm)}^n$ by $\Um$, where the $i^{\text{th}}$ factor is given the
symplectic structure
which is $2r_i$ times the usual Fubini--Study form.
In the next subsection we will describe the corresponding {\em analytic}
quotient
of $({(\bbCPm)}^n$ by $\GLm$
in the sense indicated above. In particular, when $\br$ is integral
then $M_{\br}$ will have a canonical structure of a complex projective
variety.

\subsubsection{Weighted semistable configurations on $\bbCPm$}
In this subsection we describe the semistable configurations on
${(\bbCPm)}^n$ equipped with the $\br$--dependent symplectic structure
just described. Set $[n]=\{1,2,\ldots,n\}$.

\begin{definition}
A configuration of $n$ points on $\bbCPm$ is a map $f$
from $[n]$ to $\bbCPm$.
\end{definition}

Let $\nu_{\br}$ be the measure on $[n]$ that assigns mass
$r_i$ to the point $i$. Also, recall that we have defined
$\rho=\sum r_i$.

The proof of the following theorem is left to the reader. 
For the case of integral weights it is one of the standard results in
Geometric Invariant Theory. For example, when all weights are
$\frac{1}{n+1}$ (so $\rho =1$),
it is proved in \cite[Definition 3.7/Proposition 3.4]{MumfordFogartyKirwan}.

\begin{theorem}\label{semistabilityinequalities}
A configuration $f$ on $\mathbb{CP}^m$ is semistable if and only if

$$f_*\nu_{\br}(L) \leq \frac{ \dim{L}+1}{m+1}\, \rho $$
for any linear subspace $L$ of $\bbCPm$
\end{theorem}

\begin{remark}
In this inequality, the left side is the mass of the closed subset
$L\subset \bbCPm$ for the push--forward measure $f_*\nu_{\br}$.
Intuitively, these semistability inequalities say that not too many points
can coincide, not too many can lie on a line, not too many on a plane etc.
\end{remark}

The configuration $f$ is said to be ``in general position"
if no two points coincide, at most two points lie on a line, at most
three on a plane, $\ldots$, at most $k+1$ lie in a projective subspace of
dimension $k$. The set of configurations in general position is a nonempty
Zariski--dense open subset of ${(\bbCPm)}^n$. In particular,
such configurations
exist. Thus the result that $\br$ satisfies strong triangle inequalities
implies $\Mr$ is nonempty is an immediate consequence of the following

\begin{proposition}
Suppose that $f$ is in general position. Then $f$ is weighted semistable
if and only if $\br$ satisfies the generalized triangle inequalities.
\end{proposition}
\begin{proof}
Since any subset of $k$ points in general position spans a
projective subspace of projective dimension $k-1$
it is an immediate consequence of Theorem \ref{semistabilityinequalities}
that $f$ is semistable if and only if for all $I\subset \{1,2,\cdots,n\}$
$$\sum_{i\in I}r_i \leq \frac{|I|}{m+1}\, \rho.$$
Clearly, the resulting system of inequalities contains,
and is implied by, the subset in which $|I| = 1$:

$$r_i \leq \frac{1}{m+1\,} \rho.$$

We have already noted (equation (\ref{rhotriangle})) that this system
is equivalent to the system of
strong triangle inequalities.
\end{proof}

Since configurations in general position always exist, we have
the missing implication in Theorem \ref{nonempty}.

\begin{corollary}
$M_{\br}$ is nonempty if $\br$ satisfies the strong triangle
inequalities.
\end{corollary}

\section{Smoothness of the Moduli Spaces}\label{walls}
In this section we give a sufficient condition in terms of
$\br$ for the space $\Mr$ to be smooth.

\subsection{Decomposable Polygons}
For $m=1$, it was shown
in \cite{KapovichMillson2} that $\Mr$ will have singularities
if, and only if,
the index set $\{1,\ldots,n\}$ can be
partitioned into proper subsets $I,J$ so that
\begin{equation}\label{partition}
   \sum_{i\in I} r_i = \sum_{j\in J}r_j.
\end{equation}
Then  there exists a polygon (in Euclidean space), with
the given side lengths $\br$, that is contained in a line
segment; such a polygon was called
\emph{degenerate}. It was further proved in
\cite{KapovichMillson2}
that a polygon is a singular point of $\Mr$ if and only if it is
degenerate.
We need analogs of (\ref{partition}) and of the notion of
degenerate polygon for the case $m\ge 2$.

\begin{definition}\label{HIJk}
  For $1\le k\le m$ and $I\cup J$ a proper partition of
  $\{1,\ldots,n\}$, define the hyperplane
\[
  H_{I,J,k}=\{ \br\in\bbR_+^n \mid
   k\sum_{i\in I} r_i = (m-k+1) \sum_{j\in J} r_j\}.
\]
(Note that this reduces to (\ref{partition}) when $m=1$).
The \emph{wall} corresponding to this hyperplane is
the intersection
\[ W_{I,J,k}= H_{I,J,k}\cap C(n,m+1)\]
(cf. Definition \ref{hypersimplex}).
\end{definition}

\begin{notation}\label{rhoI}
We will write
$I=\{i_1,\ldots,i_p\},
J=\{j_1,\ldots,j_k\}$, $p+q=n$, and take $I$ and $J$ to be
ordered, $i_1<i_2<\ldots$, $j_1<j_2<\ldots$. Set $\br_I=(r_{i_1},
\ldots,r_{i_p})$, and likewise for $J$. Let
$\rho_I = \sum_{i\in I} r_i$ (similarly for $\rho_J$),
and define $\vL_I=\rho_I/(m-k+1), \vL_J=\rho_J/k$,
by analogy with $\vL=\rho/(m+1)$.
\end{notation}

\begin{lemma}\label{subtriangleinequalities}
Suppose that $\br \in W_{I,J,k}$. Then  $\br_I$ (resp. $\br_J$) satisfies the
strong triangle inequalities with weight $m-k$ (resp $k-1$).
Explicitly:
\begin{align*}
r_i \leq \  & \frac{1}{m-k+1} \  \rho_I, \ \text{for all} \ i \in I \\
r_j \leq \  & \frac{1}{k} \ \rho_J, \ \text{for all} \ j \in J.
\end{align*}
\end{lemma}

\begin{proof}
We show that $\br_I$ satisfies the strong triangle inequalities with
weight $m-k$. According to Definition \ref{HIJk}, $k\rho_I=
(m-k+1)\rho_J$. Obviously, $\rho_I+\rho_J=\sum r_i=\rho$. Solving these
two equations for $\rho_I$, we get
$(m+1)\rho_I=(m-k+1)\rho$, or
$$\frac{\rho_I}{m-k+1}=\frac{\rho}{m+1}.$$
Since $\br$ already satisfies the strong triangle inequalities of
weight $m$, equation (\ref{rhotriangle}) shows the right side to be
greater than $r_i$. This gives the desired inequality for $\rho_I$; the
proof for $\rho_J$ is similar.
\end{proof}

\begin{lemma}\label{prodpoly}
If $\br \in W_{I,J,k}$, then $\vL_I=\vL_J=\vL$.
\end{lemma}

\begin{proof}
   Since $\br\in W_{I,J,k}$, we have $k\rho_I=(m-k+1)\rho_J$,
   which implies $\vL_I=\vL_J$. Furthermore,
$$k\rho=k\rho_I+k\rho_J=(m-k+1)\rho_J+k\rho_J=(m+1)\rho_J,$$
whence $\vL=\rho/(m+1)=\rho_J/k=\vL_J$.
\end{proof}

We will see that if $\br$ does not lie on a wall, then
$\Mr$ is smooth. To this end, we need the
analog of degenerate polygon. It is the ``decomposable polygon'',
in which the edges indexed by $I$ and $J$ act in orthogonal
subspaces.

Let $\br\in W_{I,J,k}$.  Choose an
orthogonal decomposition $\bbCm=V_1\oplus V_2$ with dim$\,V_1=m-k+1$
and dim$\,V_2=k$. Let $\calH_i$ denote the
set of Hermitean endomorphisms
of $V_i$. We have inclusions $\alpha_i:\calH_i\to\calHm$ given
by
$$\alpha_i(X) = \iota_{V_i}\circ X \circ \pi_{V_i}$$
where $\iota_{V_i}$ is the inclusion of $V_i$ into $V$ and
$\pi_{V_i}$ is the orthogonal projection from $V$ to $V_i$.
Note that if $w \in V_i$, then $w \otimes w^*$ is in the image
of $\alpha_i$.

We wish to define a map
\[ \iota_{I,J,V_1,V_2}: \Mtw_{\br_I}(\calH_1)\times \Mtw_{\br_J}
   (\calH_2)\to\Nrtw.
\]

Let $\sigma$ be
the $p,q$ shuffle permutation given by

\begin{align*}
\sigma(k) & =  i_k, \ 1 \leq k \leq p \\
\sigma(p+l) & =  j_l,\  1 \leq l \leq q.
\end{align*}
Choose polygons
$$\be^{(1)} = (e_1^{(1)},\cdots, e_p^{(1)}) \in
\Mtw_{\br_I}(\calH_1),\
\be^{(2)} =
(e_1^{(2)},\cdots, e_q^{(2)}) \in \Mtw_{\br_I}(\calH_2).$$
Such polygons exist, since $\br_I$ and $\br_J$ satisfy the
strong triangle inequalities (Lemma \ref{subtriangleinequalities}).

Then we define $\overline{\be} := \iota_{I,J,V_1,V_2}(\be^{(1)}, \be^{(2)})$ by
\begin{align*}
\overline{e}_{i_k} =\  & \alpha_1(e_k^{(1)}),\  1 \leq k \leq p \\
\overline{e}_{j_l} =\  & \alpha_2(e_l^{(2)}), \  1 \leq l \leq q.
\end{align*}
We then have $\overline{\be}=\overline{\be}_I\oplus\overline{\be}_J$,
and it follows from
Lemma \ref{prodpoly} that $\overline{\be}\in\Mrtw$. Therefore the image of
$\iota_{I,J,V_1,V_2}$ lies in $\Mrtw$, i.~e.~, consists of closed
polygons.

\begin{definition} We say that $\be\in\Mrtw$ is \emph{decomposable}
if it lies in the image of the map $\iota_{I,J,V_1,V_2}$ for some
choice of $I,J,V_1,V_2$ as above.
\end{definition}

\begin{proposition} $\Mrtw$ contains a decomposable polygon if, and only if,
$\br$ lies on a wall.
\end{proposition}

\begin{proof}
Suppose that $\be$ is decomposable. We use the notation above.
Since $\be_I$ is closed, $\sum_I e_i = \Lambda_I \id_1$,
where $\id_1$ is the identity in $End(V_1)$. Similarly
$\sum_J e_j = \Lambda_J \id_2$. Because $\overline{\be}$ is closed,
we have
$\sum_{k=1}^n \overline{e}_k = \Lambda \id$ ($\id$ is the identity
on $\bbCm$). Clearly this last sum is also of block form
$\Lambda_I \id_1 \oplus \Lambda_J\id_2$ in
$End(\bbCm) = End(V_1 \oplus V_2)$.
Hence $\Lambda_I = \Lambda_J = \Lambda$. This implies
$\br \in W_{I,J,k}$.

Conversely, suppose that $\br \in W_{I,J,k}$. By Lemma
\ref{subtriangleinequalities} and Theorem \ref{nonempty}
there exist a closed $p$-gon $\be^{(1)} \in \calH_1$ with side-lengths
$\br_I$ and a closed $q$-gon $\be^{(2)} \in \calH_2$  with
side lengths $\br_J$. Then $\iota_{I,J,V_1,V_2}(\be^{(1)}, \be^{(2)}) $
is a decomposable polygon with side-lengths $\br$.
\end{proof}

The next few lemmas are well known, but we will state and prove them
because
some of them will play a critical role later. As above, $\calH_j$
denotes the Hermitean endomorphisms of a subspace $V_j\subset \bbCm$ and
$\alpha_j:\calH_j \to \calHm$ is the natural inclusion.

For $\be\in\Nrtw$,
let $Z(e_i)$ denote the centralizer of $e_i$ in $\calHmz$,
and let $Z(\be)=\bigcap_i Z(e_i)$.

\begin{lemma}\label{basiclemma}
Suppose that
$X\in Z(e_i)$, and
let $\bbCm= \oplus_{j=1}^l V_j$ be the eigenspace decomposition of $X$. Then
there exists $j$ such that $e_i \in \alpha_j(\calH_j)$.
\end{lemma}

\begin{proof}
If $e_i=r_i\wtens{i}$,
we have $\bbC w_i={\rm ker}\,(e_i-r_i\id)$. Since $X$ and $e_i$
commute, $w_i$ is also an eigenvector of $X$. Hence $w_i\in V_{j}$
for some $j$ and $r_i\wtens{i} \in \alpha_j(\mathcal{H}(V_j))$.
\end{proof}

\begin{lemma}\label{imageofderivative}
Let $\be\in\Nrtw$. Then $Im \ T_{\be}\mur \subset \mathcal{H}_{m+1}^0.$
\end{lemma}

\begin{proof} If $w(t)$ is a smooth curve in
$\bbCm$, with $\|w(t)\|\equiv1$, then $\Tr w(t)\otimes w(t)^*\equiv1$
implies $\Tr\frac{d}{dt}(w(t)\otimes w(t)^*)\equiv0$.
Hence the derivative of the momentum map $\mur:\be\mapsto
\sum_i r_i \wtens{i}$ maps
into $\calHmz$.
\end{proof}

\begin{lemma}\label{Tmuperp}
Let ${}^\perp$ denote orthogonal complement in $\calHmz$. Then,
again for $\be\in\Nrtw$,
$$(Im \ T_\be\mu)^{\perp} = Z(\be).$$
\end{lemma}

\begin{proof}
Indeed,
the derivative $T_e\mur$ will be onto if, and only if,
\[
T_{e_1}(\orb_{r_1})+\cdots +T_{e_n}(\orb_{r_n})=\calHmz.
\]
Thus $T_e\mur$ is onto if, and only if,
\[
T_{e_1}(\orb_{r_1})^\perp \cap\ldots \cap T_{e_n}(\orb_{r_n})^\perp=\{0\}.
\]
But $T_{e_i}(\orb_{r_i})^\perp=
\{[e_i,X]\mid X\in\calHmz\}^\perp=Z(e_i)$, and the lemma follows.
\end{proof}

\begin{corollary}\label{T_muOnto}
Let $\be\in\Nrtw$. Then
\[
T_{\be}\mur: T_\be(\Nrtw)\to\calHmz \hbox{\ is not onto\ }
\iff Z(\be)\ne \{0\}.
\]
\end{corollary}

The next lemma relates surjectivity of the moment map to indecomposability.

\begin{lemma}\label{NotOntoIsDecomp}
Now suppose that $\be\in\Mrtw$. Then
$$T_{\be}\mur \text{ is onto} \iff \be \ \text{is not
decomposable}.$$
\end{lemma}
\begin{proof}
Suppose that $T_{\be}\mur$ is not onto. Choose a nonzero
$X\in Z(\be)$. Suppose that $X$ has $\ell$ distinct
eigenvalues, so that $\bbCm$
is the orthogonal sum of the corresponding eigenspaces $W_j$.
For each $e_i=r_i\wtens{i}$,
$w_i\in W_{j_i}$
for some $j_i$, by Lemma \ref{basiclemma}.
Now set $V_1=W_1+\cdots+W_{\ell-1}$, $V_2=W_\ell$. Define
$I=\{i\mid w_i\in V_1\}, J=\{j\mid w_j\in V_2\}$.
It follows that $\be$ lies in the image of the map $\iota_{I,J,V_1,V_2}$
and so is decomposable.

Now suppose that  $\be$ is decomposable. Then there exists an orthogonal
splitting $V = V_1 \oplus V_2$ and a partition
$\{1,\cdots,n\}$ = $I \cup J$
such that $\be=\be_I+\be_J$  is in the
image of the map $\iota_{I,J,V_1,V_2}$. Let $X=\id_1\oplus{\mathbf 0}$.
Then $X\in Z(\be)=(\Imag \ T_{\be} \mur)^{\perp}$.
\end{proof}

Let $\sigrtw\subset\Mrtw$
denote the set of decomposable polygons. It is invariant
under $\Um$; let $\sigr$ be the image of $\sigrtw$ in $\Mr$.

\begin{theorem}\label{singpoint}
(i) $\Mrtw-\sigrtw$ is a smooth manifold. (ii) The group $\SUm$ acts
freely on $\Mrtw-\sigrtw$, hence the quotient $\Mr-\sigr$ is a smooth
manifold.
\end{theorem}
\begin{proof}
Part (i) follows from one implication in Lemma
\ref{NotOntoIsDecomp}: if $\be$ is not decomposable,
then $\vL\id$ is a regular value of $\mur$.

For (ii), we need to check that
if $\be$ is not decomposable, then the  stabilizer of $\be$
under the action of $\Um$ is trivial. The argument in Lemma
\ref{NotOntoIsDecomp} still works,
because we deal with matrix groups. If $\kappa\be \kappa^{-1}=\be$,
we write
$\bbCm$ as sum of eigenspaces of $\kappa$, and proceed as before.
\end{proof}

\begin{corollary}\label{dimM_r}
If $\br$ does not lie on a wall, then $\Mr$ is a smooth manifold of dimension
$2m(n-m-2)$.
\end{corollary}
\begin{proof}
The smoothness of $\Mr$ follows from the theorem. To compute the
dimension of $\Mr$, we note that $\Mr$ is the symplectic quotient of
${(\bbCPm)}^n$ by the projective unitary group $PU(m+1)$. Thus we
obtain
$$ \dim_{\bbR} \Mr = 2mn - 2[(m+1)^2 - 1] = 2mn - 2m^2 - 4m.$$
\end{proof}

\subsection{The critical sidelengths of closed polygons}

In this subsection, we study the map $\barbs$ that maps a closed
polygon $\be$ to the vector $\br$ of its side lengths, and show in
Theorem \ref{wallunion} that its critical values are the union of
the walls (Definition \ref{HIJk}).

\begin{notation}
The space of not necessarily closed linkages
in $\calHm$ with rank $1$ positive semi-definite
edges and
{\em arbitrary (positive) side lengths} is denoted by
$$\Pol=\{\be\mid \be\in\Nrtw \text{\ for some\ }\br\}. $$
We will use $\mu$ to denote the restriction of the
momentum map $\mu$ for the diagonal action of $U(m+1)$ on $\calHm^n$
to $\Pol$. (Recall that $\mur$ is the moment map on linkages, closed
or not, with given side lengths $\br$).
Further, introduce the subset of {\em closed} polygons
with arbitrary side lengths,
$$\CPol=\{\be\mid \be\in\Mrtw \text{\ for some\ }\br\}\subset \Pol. $$
\end{notation}

The idea of the argument is this. If
$\be(t)=(r_1(t)w_1(t)\otimes w_1(t)^*,
\ldots)$ is a curve in $\CPol$ through $\be(0)=\be,\br(0)=\br$,
the derivative of the moment map $\mu$ has the form
$$T_\be\mu=\sum \dot r_i(0)e_i + \sum r_i(w_i\otimes \dot w_i(0)^*+
\dot w_i(0)\otimes w_i^*).$$
The first sum is a linear combination of edges, while the second sum
is in the image of the moment map $\mur$. Since
$T_\be\barbs=\dot \br(0)$, the span of the edges relates
surjectivity of $T_\be\mu$ with surjectivity of $T_\be\barbs$.

Let $\mathbf{E} \subset \calHm$ be the span of the edges
$e_i$ of the linkage $\be$.

\begin{lemma}\label{centralizerinspan}
Suppose that $\be \in CPol$. Then
 $$Z(\be) \subset \mathbf{E}.$$
\end{lemma}

\begin{proof}
Let $X \in Z(\be)$, with eigenvalues $\lambda_i,i=1,\ldots,l$, and
let $\bbCm= \oplus_{i=1}^l V_i$ be the eigenspace splitting of $\bbCm$
under $X$.
Then by Lemma \ref{basiclemma}, for each $i,1 \leq i \leq n$, there exists
$j_i$ such that $e_i \in \alpha_{j_i}(\calH_{j_i})$.
Hence, $\be$ is decomposable with respect to this splitting.
Thus there exists a permutation $\sigma$
such that
$ \be = \sigma(e^{(1)}_1,\cdots,e^{(1)}_{p_1}, \cdots, e^{(l)}_1,\cdots,
e^{(l)}_{p_l})$ with $e^{(i)}_j=r{(i)}_j w{(i)}_j\otimes (w^{(i)}_j)^*$, and
$w^{(i)}_j\in V_i$,\ $1\leq i \leq l$, $1\leq j \leq p_i$.

Since $\be$ is closed, we have
$\sum_i e_i = \vL \id$. As a consequence we have
$$\sum_{j=1}^{p_i} e^{(i)}_j = \vL \Pi_i$$
where $\Pi_i$ is the projection on $V_i$.
Thus the
$\Pi_i$ lie in the span $\mathbf{E}$. But since the splitting of
$\bbCm$ is
the eigenspace decomposition of $X$, we have
$X = \textstyle{\sum_{i=1}^l \lambda_i \Pi_i}$.
Thus $X \in \mathbf{E}$.
\end{proof}

We do not need the next corollary in what follows but have included it
for completeness. We let $\mathfrak{t}$ denote
the abelian Lie subalgebra consisting of the diagonal matrices in
$\calHm$.

\begin{corollary}
There exists $k \in U(m+1)$ such that
$$\Ad_k(Z(\be)) \subset \mathfrak{t}.$$
\end{corollary}

\begin{proof}
It suffices to prove that $Z(\be)$ is abelian. To this end let $X,Y \in
Z(\be)$. By the lemma we may write $X$ as a linear combination of the
edges of $\be$. But by definition $Y$ centralizes all the edges of $\be$.
\end{proof}

\begin{lemma}\label{spanisenough}
$$\Imag\ T_\be \mu = \Imag\ T_e \mu_r + \mathbf{E}.$$
\end{lemma}

\begin{proof}
Define an action of $\bbRnp$ on $\Pol$ by
$\mathbf{a} \cdot \be = (a_1e_1,a_2 e_2 \cdots, a_n e_n)$, where
$\mathbf{a} = (a_1,\cdots, a_n)$. It is immediate that $\Nrtw$ is a
cross-section
for this action and consequently we have
$\Pol \cong \bbRnp\times \Nrtw$.  Now let $\be \in \Pol$.
Then
$$\mu(\mathbf{a}\cdot \be) = a_1e_1 + a_2 e_2 +\cdots + a_n e_n.$$
The lemma follows upon differentiating this identity with respect to
$\be$ and the action of $\bbRnp$.
\end{proof}

The following proposition will play a critical role in our analysis of the
diagram below.

\begin{proposition}\label{ontoprop}
Suppose that $\be\in\CPol$.
Then $T_{\be}\mu$ maps onto $\calHm$.
\end{proposition}

\begin{proof}
We again use the fact that
$(\Imag \ T_{\be} \mu_{\br})^\perp =Z(\be)$. By Lemma
\ref{centralizerinspan} the
directions coming by changing the side lengths (the action of $\bbRnp$)
contain $Z(\be)$. The result now follows from Lemma \ref{spanisenough}.
\end{proof}

Let $\bs:\Pol \to\bbRn$
and $\overline{\bs} :\CPol\to\bbRn$ denote the  side-length maps,
and let ${\bf j}:\CPol \to \Pol$ and
${\bf k}:\Nrtw \to \Pol$ be the inclusions. It is immediate
(by using the action of $\bbRnp$) that $T_{\be} \bs$ maps onto $\bbRn$.

We will need the subspaces of $\Pol$ and $\CPol$ obtained
by fixing the sums of the side lengths (but not the side lengths
themselves). For $\lambda\in\bbR_+$, define
\begin{align*}
\Pol_{\lambda} &= \{ \be \in \Pol \mid \sum\limits_{i=1}^n
||e_i|| = \lambda \}\\
\CPol_{\lambda}&= \CPol\cap \Pol_{\lambda}.
\end{align*}
We observe that $T_\be \mu$ maps $T_{\be}\Pol_{\lambda}$
into $\calHm^0$, and moreover
it is an immediate consequence of Proposition \ref{ontoprop} that
if $\be$ is closed then this map is onto.

\begin{remark} \label{exactness}
In what follows we will use the fact that if $f:X \to Y$ is
a real-analytic
map of real analytic sets with $f(x) =y$ then the sequence
$$T_x(f^{-1}(y)) \to T_x(X) \to T_y(Y)$$
is exact at $T_{x}(X)$ (the second arrow is $T_x(f)$).
\end{remark}

\begin{theorem}\label{wallunion}
The set of critical values of $\barbs$ is the union of the walls.
\end{theorem}

\begin{proof}
We have seen that $\br$ lies on a wall if, and only if, $\vL_\br\id$
is a critical value of $\mu_{\br}:\Nrtw\to\calHm$. Let
$\be\in\Mrtw$.
The result will follow once we prove that $T_\be \mu_{\br}$ is onto
if, and only if, $T_{\be}\overline{\bs}$ is onto. This follows from a diagram
chase
in the following commutative diagram. Recall
that $\rho =\textstyle{ \sum_{i=1}^n r_i}$. Let $\bbRn_0$ denote
the subspace of $n$-tuples with sum $0$. It is clear that
$T_{\be}\bs$
maps $T_{\be}\Pol_{\rho}$ onto $\bbRn_0$.

\begin{center}
\begin{math}
\begin{CD}
0 @>>> T_{\be}\Mrtw @>T_{\be}{\bf i}>> T_{\be}\Nrtw @>T_{\be}\mu _{\br}>>
\calHm^0 \\
@VVV   @VVV                @VVT_{\be}{\bf k}V   @VV{\id}V \\
0 @>>> T_{\be}\CPol_{\rho} @>T_{\be}{\bf j}>> T_{\be}\Pol_{\rho}
@>T_{\be}{\mu}>> \calHm^0\\
 @VVV  @V{T_{\be} \overline{\bs} } VV  @VV{ T_{\be}{\bs}}V  @VVV\\
 0 @>>> \mathbb{R}^n_0 @>{\id}>> \mathbb{R}^n_0 @>>> 0
\end{CD}
\end{math}
\end{center}

\vspace{ .25 cm}

We will perform  the  diagram chase that proves  $T_\be\overline{\bs}$  onto
$\Longrightarrow  T_{\be}\mu _{\br}$  onto. To this end let
$y \in \calHmz$.
We will construct $x \in T_{\be}\Nrtw$ with  $T_{\be}\mu _{\br}(x) = y$.
Indeed, since $T_{\be}\mu $ is onto there exists
$z \in  T_{\be}\Pol_{\rho}$ with $T_{\be}\mu (z) = y$.
By our assumption that $T_{\be}\barbs$ is onto there exists
$w \in T_{\be}\CPol_{\rho}$ with $T_{\be}\barbs(w) =
T_{\be}{\bs}(z)$.
Then  $T_{\be}\bs(z - T_{\be}{\bf j}(w)) = 0$.

Since the next to last vertical sequence is exact at $T_{\be}\Pol$
by Remark \ref{exactness}, there exists $x \in T_{\be}\Nrtw$ with
$T_{\be}{\bf k}(x) = z - T_{\be}{\bf j}(w)$. Then
 $T_{\be}\mu _{\br}(x) = y$ as required.
 \end{proof}

This theorem is a critical first step towards finding the topologies of
the moduli spaces $\Mr$. We define a \emph{chamber} of the
polyhedral cone
$C(n,m+1)$ to be a connected component of the complement of the
union of the walls.
Then we have the following corollary of the previous theorem.

\begin{corollary}
The sidelength map $\mathbf{s}$ is a (trivial) fiber bundle over each
chamber of $C(n,m+1)$. Hence if $\br$ and $\br^{\prime}$ lie
in the same chamber, the moduli spaces $\Mr$ and $M_{\br^{\prime}}$
are diffeomorphic.
\end{corollary}
\begin{proof}
The map $\mathbf{s}$ is a proper submersion over each chamber,
hence by the
Ehresmann Fibration Theorem, \cite[pg. 84]{BrockerJanich}, it is a
fiber bundle (necessarily trivial since chambers are contractible).
\end{proof}

In fact we have to replace $\Pol_{\rho}$ by its quotient
by $SU(m+1)$ in the argument above, but the reader will check that
the critical set of the map induced by $\barbs$ on this
quotient remains the same.

\begin{remark}
It is possible to implement wall-crossing techniques in order to compute the 
topologies of
the moduli spaces $\Mr$ --see\cite{Hu} or \cite{Goldin} for
cohomology computations.
\end{remark}

\section{Bending Hamiltonians}\label{bending}
Kapovich and Millson (\cite{KapovichMillson2}) studied an integrable
Hamiltonian system on  $\Mrtw$ in the case $m=1$ and $e_2=-e_1$,
i.e. in $\calH_2^0\equiv\sutwo$,
which is isomorphic to Euclidean space $\bbE^3$.
In \S1 we introduced the diagonals
$A_0=e_1$ and $A_i=e_1+\cdots+e_{i+1},i=1,\ldots,n-2 $.
In this case, a closed polygon has $A_{n-1}=\vL\id$.
It was shown that the functions $f_i(\be)=\|A_i\|$
Poisson commute; $\|A_i\|$ is the positive eigenvalue of $A_i$.
The diagonal $A_i$ divides the polygon
into two ``flaps'', and the flow generated by $f_i$ is $2\pi$-periodic,
consisting of a rigid rotation of one flap about the diagonal.

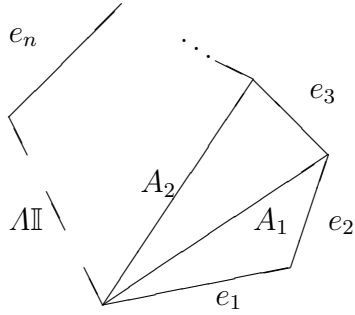
\begin{figure}[h]
\begin{center}
\setlength{\unitlength}{5mm}
\begin{picture}(12,12)(-4,-1)
\linethickness{1.5pt}
\put(0,0){\line(5,1){5}}
\put(5,1){\line(1,3){1}}
\put(6,4){\line(-1,1){2}}
\put(4,6){\line(-2,1){1}}
\put(0,0){\line(-1,2){.5}}
\put(-1,2){\line(-1,2){.5}}
\put(-2,4){\line(-1,2){.5}}
\put(-2.5,5){\line(1,1){3}}
\put(0,0){\line(3,2){6}}
\put(0,0){\line(2,3){4}}
\put(3,0){$e_1$}
\put(6,2){$e_2$}
\put(2,6.5){$\ddots$}
\put(5.5,5.5){$e_3$}
\put(-2.5,7){$e_n$}
\put(-2.5,2){$\vL\mathbb{I}$}
\put(4,2){$A_1$}
\put(1,3){$A_2$}
\end{picture}
\end{center}
\caption{A polygon in $\mathfrak{u}(m+1)$}
\end{figure}

The analogs of ``bending Hamiltonians'' for $m>1$
are again the eigenvalues of the diagonals. Now, however, $A_{n-1}=
\vL\id$, indicated by a dashed line in Figure~1; it would be
absent in ${\mathfrak su}(m+1)$.

\begin{nnotation}\label{eivals}
The eigenvalues of $A_i$ are denoted by $\laij$ in decreasing
order, $\lambda_{i1}\ge\ldots\ge\lambda_{i,m+1}$.
\end{nnotation}

We note that $A_{n-2}=\vL\id-e_n$, has eigenvalues $\vL$
(multiplicity $m$)
and $\vL-r_n$, and those are fixed. Thus only the $\laij$ for
$1\le i\le n-3$ are of possible interest. Furthermore, it
will be seen in \S\ref{WA} that off submanifolds
of $\Mrtw$ of lower dimension, the nontrivial $\laij$
(those not identically $0$ or $\vL$) are simple. In that
case, they will be smooth functions of $\be$, which is
assumed throughout the present section.

\subsection{Bending Flows}
We want to calculate the Hamiltonian vector fields and flows
generated by the $\laij$. By analogy with the case of $\mathbb{E}^3$,
we call them ``bending flows''.

On a product of orbits, the Poisson bracket is the sum
of the orbit brackets, and the next formula is evident from
(\ref{Ham-equn}):
\begin{proposition}\label{kirillovsum}
Suppose $f:\Mrtw\to\bbC$ is smooth and depends
only on $e_1,\dots,e_{i+1}$. Then the Hamiltonian system generated by
$f$ is
\begin{align}\label{Hamequs}
\dot e_k =
\begin{cases}
 [\nabla_k f(e_1,\dots,e_{i+1}),e_k],&\text{if $1\le k\le i+1$,}\\
 0,                                   &\text{if $i+1< k\le n$,}
\end{cases}
\end{align}
where $\nabla_k$ denotes gradient with respect to $e_k$, all other
$e_j$ being held fixed.
\end{proposition}

To solve these equations when $f=\laij$, we recall a standard lemma from
perturbation theory.
\begin{lemma}\label{perturb}
Let $\lambda$ be an isolated eigenvalue of $A\in\calHm$, with
unit eigenvector $u$. Then $\nabla\lambda(A)=\ei u\otimes u^*$.
\end{lemma}
\begin{proof}
For $A'$ sufficiently close to $A$, the eigenvalue $\lambda(A')$ and
(with proper choice of phase)
normalized eigenvector $u(A')$ vary analytically in a neighborhood
of $\lambda,u$. Take a curve
$A(s)u(s)=\lambda(s)u(s)$, and take the inner product with the unit
length $u(s)$ to get $\lambda(s)=(A(s)u(s),u(s))$. Differentiate
and set $s=0$, and use $(Au,\dot u(0))+(A\dot u(0),u)=\lambda
((u,\dot u(0))+(\dot u(0),u))=0$, resulting in
\[
\dot\lambda(0)=(\dot A(0)u,u)=\Imtr(\ei u\otimes u^*\dot A(0)),
\]
as was to be shown.
\end{proof}

Write $E_j(A)$ for the spectral projection onto the $\lambda_j$
eigenspace of $A$; the lemma thus states that $\nabla \lambda_j(A)=
\ei E_j(A)$. We compute the Hamiltonian flows generated by the
functions $\laij:\be\mapsto\lambda_j(A_i)$ on
$\Nrtw$.  These functions are invariant under the $\Um$ action,
and hence descend to the symplectic quotient $\Mr$. The
$\laij$ are smooth when they are simple eigenvalues.

\begin{proposition}\label{phi-ij}
For $i=1,\ldots,n-3$ and $j=1,\ldots,m+1$,
$\laij$ is the Hamiltonian for
the system
\begin{align} \label{phi-ij-equ}
\dot e_k =
\begin{cases}
 \ei[E_j(e_1+\dots+e_{i+1}),e_k],   & \text{if $1\le k\le i+1$,}\\
 0,                                   & \text{if $i+1< k\le n$.}
\end{cases}
\end{align}
The Hamiltonian flow $\phiij^t(\be)=\be(t)$ is given by
\begin{align} \label{phi-ij-flow}
 e_k(t) =
\begin{cases}
\big(\Ad \exp(\ei t E_j(A_i))\big)(e_k),   &\text{if $1\le k\le i+1$,}\\
 e_k,                                   &\text{if $i+1< k\le n$.}
\end{cases}
\end{align}
\end{proposition}
\begin{proof}
To obtain the system (\ref{phi-ij-equ}) we wish to apply
Proposition \ref{kirillovsum}. It is necessary to relate
the partial gradients $\nabla_k\laij$ to the full gradient,
$\nabla\laij=\ei E_j(A_i)$. According to Lemma \ref{perturb},
the former are found by computing
\[
\dot A_i(s)=(e_1+\dots+e_k(s)+\dots+e_{i+1})\dot{}=\dot e_k(s),
\]
but because $\dot A_i(0)$ is tangent to $\orb_{r_k}$, this only
determines $\nabla_k\laij$ up to a vector normal
to the orbit:
\[
\nabla_k\laij(A_i)=\ei E_j(A_i)+\xi_k,\quad [\xi_k,e_k]=0.
\]
Then $[\nabla_k\laij,e_k]=\ei [E_j(A_i),e_k]$, and (\ref{phi-ij-equ}) follows.


Next, add the equations (\ref{phi-ij-equ}) for $1\le k\le i+1$ to find
\[\dot A_i(t)=
\ei [E_j(A_i(t)),A_i(t)].
\]
Since $A_i$ commutes with its own spectral
projections, we get $\dot A_i(t)=0$ and $A_i(t)=A_i$. With constant
$A_i$, the solution of (\ref{phi-ij-equ}) is immediate.
\end{proof}
\begin{corollary}\label{twopi}
The flows $\phiij$ have period $2\pi$ in $t$.
\end{corollary}
\begin{proof}
If $P$ is a projection, then $P^2=P$. Consequently, $\exp(\ei tP)=
\id +(\exp(\ei t) -1)P$, which has period $2\pi$.
\end{proof}

\subsection{Involutivity}
It is not \emph{a priori} clear from the formulas for $\phiij$ that
these flows commute. This is a short calculation; we again work only
with simple eigenvalues of the $A_i$ on $\Nrtw$, and the flows
will also commute on $\Mr$.

\begin{proposition}
$\{\laij,\lakl\}=0$ for $1\le i,k\le n-3$ and $1\le j,\ell\le m+1$.
\end{proposition}
\begin{proof}

By Proposition \ref{kirillovsum} and the proof of Proposition \ref{phi-ij},
\[
\{\laij,\lakl\}(\be)=\sum_{s=1}^{i+1}\Imtr\bigl(e_s
[\ei E_j(A_i)+\xi_s,\ei E_\ell(A_k)+\eta_s]\bigr),
\]
where again $\xi_s,\eta_s$ commute with $e_s$. The ad-invariance
of the trace form produces $[\xi,e_s]$ and $[\eta,e_s]$, which are zero. This leaves
\begin{align*}
\{\laij,\lakl\}(\be)&=-\sum_{s=1}^{i+1}
\Imtr\bigl([e_s,E_j(A_i)]E_\ell(A_k)\bigr)\\
                    &=-\Imtr\bigl([A_i,E_j(A_i)]E_\ell(A_k)\bigr)\\
                    &=0.
\end{align*}
\end{proof}

\begin{remark}
The proof works more generally, if instead of $A_i$ and $A_k$ one has
$\sum_I e_i$ and $\sum_J e_j$, with $I\subset J$. Thus, for example,
the eigenvalues of $e_2+e_3$ and $e_1+\cdots+e_5$ are in involution.
On the other hand, if $\lambda,\mu$ are eigenvalues of $e_1+e_2$ and $e_2+e_3$,
respectively, then
\[
\{\lambda,\mu\}(\be)=-\Imtr\bigl(e_2[E_\lambda(e_1+e_2),
E_\mu(e_2+e_3)]\bigr),
\]
which need not be zero.
See \cite{KapovichMillson3} for more information.
\end{remark}

\section{A Complete Set of Bending Flows}

The eigenvalues $\laij(\be)$ have been shown to Poisson commute,
and to generate $2\pi$-periodic flows. If there were
$\frac{1}{2}\dim \Mr$ eigenvalues and if they were smooth, they
would constitute a set of action variables on $\Mr$.
Smoothness everywhere cannot be achieved, but there
are $\frac{1}{2}\dim\Mr$
that are smooth and functionally independent on a dense open
submanifold of $\Mr$. This section presents the proof.

\subsection{The Weinstein-Aronszajn Formula}\label{WA}

The diagonal $A_i$ is a rank-one perturbation of $A_{i-1}$, and
because of this, the eigenvalues $\laij$ and $\lambda_{i-1,j}$ are
related in a special way. This connection is the simplest instance
of the Weinstein-Aronszajn formula
\cite[Ch.4, \S6]{Kato}. We describe the
formula and two consequences that will be used later.

Let $A$ be an $(m+1)\times(m+1)$ Hermitean matrix with eigenvalues
$\lambda_1,\dots,\lambda_{m+1}$ and let $u_1,\dots,u_{m+1}$
be corresponding orthonormal eigenvectors. (If an eigenvalue
has multiplicity $>1$, which is now permitted, the choice of
its eigenvectors is irrelevant). Let $w\in\bbCm$ be a unit vector, and let
$r\in\bbR$. Set $L=A+r\wtens{{}}$, and call its eigenvalues
$\nu_1,\dots,\nu_{m+1}$. Finally, define $\alpha_1,\dots,
\alpha_{m+1}\in\bbC$ by $w=\sum_{j=1}^{m+1}\alpha_j u_j$.

\begin{proposition}\label{detratio}
  \begin{equation}\label{detratio-equ}
  \frac{\det(z\id-L)}{\det(z\id-A)}=1-r\sum_{j=1}^{m+1}
  \frac{|\alpha_j|^2}{z-\lambda_j}.
  \end{equation}.
\end{proposition}
  \begin{proof}
Write $R_z=(z\id-A)^{-1}$ for the resolvent of $A$. The left side
of (\ref{detratio-equ}) is
    \begin{align*}
&=\det\big((z\id-A)^{-1}(z\id-A-r\wtens{{}})\big)\\
&=\det(\id -R_z(r\wtens{{}}))\\
&=\det(\id-r(R_zw)\otimes w^*).
    \end{align*}
Now, $\det(\zeta\id-r(R_zw)\otimes w^*)$ is the characteristic
polynomial of a rank-one matrix, and so has an $m$-fold root at
$\zeta=0$ and a simple root at $\zeta=r(R_zw,w)$. Setting $\zeta=1$
we get
   \begin{equation}\label{resolvent}
   \det(\id-r(R_zw)\otimes w^*)=1-r(R_zw,w).
   \end{equation}
The lemma now follows by expanding $w$ in (\ref{resolvent}) in the basis $u_j$.
   \end{proof}

It is convenient to write (\ref{detratio}) more explicitly:
   \begin{equation}\label{prodratio}
     \frac{(z-\nu_1)\dots(z-\nu_{m+1})}
     {(z-\lambda_1)\dots(z-\lambda_{m+1})}=
     1-r\sum_{j=1}^{m+1}\frac{|\alpha_j|^2}{z-\lambda_j}.
   \end{equation}
\begin{corollary}\label{rationalfunctions}
The $|\alpha_j|^2$ are rational functions of $\nu_k,\lambda_\ell,
1\le k,\ell\le m+1$.

\end{corollary}

Finally, we show that the eigenvalues of $A$ and $L$ interlace.
This will play a basic role below.

\begin{proposition}\label{interlace}
  If $r>0$, then $\nu_1\ge\lambda_1\ge\nu_2\dots\ge\nu_{m+1}\ge
\lambda_{m+1}$.
  If $r<0$, we have $\lambda_1\ge\nu_1\dots$ instead.
\end{proposition}
   \begin{proof}
Suppose $r>0$.
It suffices to prove the proposition for a dense set of $w$, so that
we may assume $|\alpha_j|^2>0$ for all $j$. Let $R(z)$ be the
rational function on the right side of (\ref{prodratio}). Since
$\lim_{z\to\infty}R(z)=1$ and $\lim_{z\downarrow\lambda_1}=-\infty$,
$R$ has a zero in $(\lambda_1,\infty)$. Likewise, because
$\lim_{z\uparrow \lambda_j}=+\infty$ and $\lim_{z\downarrow\lambda_{j+1}}
=-\infty$, $R$ has a zero in $(\lambda_{j+1},\lambda_j)$. This provides
$m+1$ zeros of $R$, which must coincide with the zeros $\nu_j$ of the
left side of (\ref{prodratio}).
   \end{proof}

\subsection{Gel'fand-Tsetlin Patterns}\label{actions}
Let $\be\in\Mrtw$. We will  arrange the eigenvalues of $A_0=e_1,
A_1,\ldots,A_{n-1}$ in a triangle with vertex at the bottom.
The eigenvalues of $A_k$ are written in row $k$ of the triangle,
along with some space-filling zeros.
For $0\le k\le m$, the rank of $A_k$ is at most $k+1$, so zero must
be at least an $(m-k)$-fold eigenvalue of $A_k$. Those zeros are not
recorded. When $k>m$, there are $m+1$ eigenvalues, potentially
nonzero; these are recorded \emph{along with} $k-m$ \emph{zeros}.
Figure~2 shows the case $m=2, n=6$. Note that entries of
successive rows are offset to reflect the interlacing property
deduced in Proposition~\ref{interlace}. This diagram is called
a \emph{Gel'fand-Tsetlin pattern}, or \emph{GTs pattern} for short.
It is denoted by $\Gamma(\be)$.
The extra zeros will be explained in \S 11, see Remark
\ref{zerogts}.
\[
\begin{array}{ccccccccccc}
\vL&{}&\vL&{}&\vL&{}&0&{}&0&{}&0\\
{}&d_1&{}&d_2&{}&d_3&{}&0&{}&0&{}\\
{}&{}&c_1&{}&c_2&{}&c_3&{}&0&{}&{}\\
{}&{}&{}&b_1&{}&b_2&{}&b_3&{}&{}&{}\\
{}&{}&{}&{}&a_1&{}&a_2&{}&{}&{}&{}\\
{}&{}&{}&{}&{}&r_1&{}&{}&{}&{}&{}
\end{array}
\]
\centerline{Figure 2}

Since $\be\in\Mrtw$, there are  additional restrictions on the
entries of $\Gamma(\be)$. Row $n-1$ must consist of $m+1$
$\vL$'s (because $e_1+\dots+e_n=\vL\id$) and $(n-m-1)$ zeros.
The interlacing property forces the first $m$ entries of row
$n-2$ to be $\vL$, so in Figure~2, $d_1=d_2=\vL$. Likewise,
$c_1=\vL$. It becomes apparent that the extra zeros remind one
that (for example) the eigenvalues $d_3=\lambda_{4,3}$ and
$c_3=\lambda_{3,3}$ must be non-negative.

Moreover,
\begin{equation}\label{rowsum}
\Tr A_k=\Tr(e_1+\dots+e_{k+1})=r_1+\dots+r_{k+1},
\end{equation}
which is a linear constraint on the rows of $\Gamma(\be)$.
In Figure~2, that leaves $c_2,b_1,b_2,a_1$ as potentially
independent commuting Hamiltonians, and indeed $\dim_\bbR\Mr=8$
in this case.

We summarize this discussion.
\begin{definition}\label{bigP}
Let $m,n,\br$ be fixed. We write $\bP$ for the convex polytope of
GTs patterns satisfying the following conditions.
\begin{enumerate}
\item There are $n$ rows numbered $0,\dots,n-1$ (starting at the bottom);
\item Row $n-1$ consists of $m+1$ $\vL$'s and $n-m-1$ zeros;
\item The sum of the entries of row $k$ is $\sum_{i=0}^k r_{i+1}$.
\item The interlacing property $\laij\ge\lambda_{i-1,j}\ge\lambda_{i,j+1}$
holds.
\end{enumerate}
\end{definition}

\begin{proposition}\label{bigPdim}
$\dim\bP=(n-m-2)m=\frac{1}{2}\dim_\bbR\Mr$.
\end{proposition}
\begin{proof}
There are two cases: (1) $n\ge2(m+1)$ and (2) $n\le2m+1$. The difference
comes from the position of row $m$, corresponding to the eigenvalues
of $A_m=e_1+\dots+e_{m+1}$. Generically, this matrix will have full
rank. In case (2), some of its eigenvalues are forced, by interlacing,
to be $\vL$. In case (1), all the automatic $\vL$'s have been ``exhausted''.
(Figure~2 falls into the latter category). Let us sketch the counting.

Case (1): Unconstrained $\laij$ can appear in rows $i=1,\dots,n-3$. Break
this index set into three parts: $S_1=\{1,\dots,m\}$, $S_2=\{m+1,\dots,
n-m-2\}$, $S_3=\{n-m-1,\dots,n-3\}$.
If $n=2(m+1)$ (as in Figure~2), then $S_2=\emptyset$.
The numbers of unconstrained $\laij$
for the corresponding $A_k$ are
\begin{itemize}
\item In $S_1$, $1,\dots,m$;
\item in $S_2$, $m,\dots,m$;
\item in $S_3$, $m-1,\dots,1$.
\end{itemize}
Adding, we obtain
\[
\frac{m(m+1)}{2}+(n-(2(m+1))m+\frac{m(m-1)}{2}=(n-m-2)m.
\]
Case (2): We set $S_1=\{1,\dots,n-m-2\}$, $S_2=\{n-m-1,\dots,m\}$,
$S_3=\{m+1,\dots,n-3\}$ (if $m=1,2$, then $S_3=\emptyset$).
The numbers of unconstrained $\laij$ are:
\begin{itemize}
\item In $S_1$, $1,\dots,n-m-2$;
\item in $S_2$, $n-m-2,\dots,n-m-2$;
\item in $S_3$, $n-m-3,\dots,1$.
\end{itemize}
Now add.
\end{proof}

\subsection{Constructing a polygon with given GTs pattern}\label{inverse}
In the last section, we saw that $\Gamma(\Mr)\subset\bP$. We now prove
the converse.
\begin{theorem}\label{inversetheorem}
(i) $\Gamma(\Mr)=\bP$. (ii) There
are $\frac{1}{2}\dim\Mr$ functionally independent $\laij$'s.
\end{theorem}

\begin{proof}
Let $\calSm\subset\calHm$ denote the space of real symmetric matrices,
and let $\Mrtw(\calSm)$ be the set of polygons in $\Mrtw$ with each
$e_i\in\calSm$. The obvious inclusion $\Mrtw(\calSm)\hookrightarrow
\Mrtw$ is the analog of the inclusion $\Mrtw(\bbR^2)\hookrightarrow
\Mrtw(\bbR^3)$ used in \cite{KapovichMillson2}. We will see later
that elements
of $\calSm(\Mrtw)$ can be thought of as ``unbent'' polygons;
these will be
important in our proof of the involutivity of the angle variables in the
next section.
We now show that
\begin{equation}\label{duister}
\Gamma(\Mrtw(\calSm))=\bP.
\end{equation}

Since $\Gamma:\Mrtw(\calSm)\to\bP$ is continuous
(though not differentiable),
the image of $\Gamma$ is closed, and it suffices to prove that the
image of $\Gamma$ contains
the interior $\bP^o$ of $\bP$. Thus, choose a GTs pattern $\gamma$ in
which all unconstrained inequalities are strict; we are to find
$\be\in\Gamma(\Mrtw(\calSm))$ such that
$\Gamma(\be)=\gamma$.

Set $A_0=r_1 w_1\otimes w_1^*$, where $w_1$ is an arbitrary real unit
vector. Assuming that a real symmetric $A_{k-1}$ with a given
spectrum has been found, we want $w_{k+1}\in\bbR^{m+1}$ so that
\begin{equation}\label{nextAk}
A_k=A_{k-1}+r_{k+1}\wtens{k+1}
\end{equation}
has the required next spectrum.

We carry out the induction step for
Case (1), in the terminology of Proposition \ref{bigPdim}.
First, let $k\in S_1$. Thus
\[A_{k-1}=\sum_{j=1}^k r_j\wtens{j};\]
it has spectrum
$\{\lambda_1,\dots,\lambda_k,0,\dots,0\}$ with $\lambda_1>\dots>\lambda_k>0$,
and $\sum_{i=1}^k\lambda_i=\sum_{i=1}^k r_i$.
We are further given
$\nu_i$ with
\[
\nu_1>\lambda_1>\nu_2>\dots>\lambda_k>v_{k+1}>0,
\]
and $\sum_{i=1}^{k+1}\nu_i=\sum_{i=1}^{k+1}r_i$.

Let $u_1,\dots,u_k, u$ be normalized (real) eigenvectors of $A_{k-1}$
corresponding to $\lambda_1,\ldots,\lambda_k,0$, and seek $w_{k+1}$
in the form
\[ w_{k+1}=\sum_{j=1}^k\alpha_j u_j + \alpha u\]
with $\alpha_i,\alpha$ real.

Now solve for $|\alpha_j|^2,1\le j\le k$ and $|\alpha|^2$ in equation
(\ref{prodratio}), which takes the special form
\[
\frac{(z-\nu_1)\dots(z-\nu_{k+1})z^{m-k}}
     {(z-\lambda_1)\dots(z-\lambda_{k})z^{m-k+1}}=
     1-r_{k+1}\biggl(\sum_{j=1}^{k}\frac{|\alpha_j|^2}{z-\lambda_j}
          \,+\,\frac{|\alpha|^2}{z}\biggr).
\]
Clearly one can take $\alpha_i,\alpha$ real.
Taking traces in equation (\ref{nextAk}), we get
\[
\sum_{j=1}^{k+1}r_j=\sum_{j=1}^{k+1}\nu_j=\sum_{j=1}^kr_j+r_{k+1}
\|w_{k+1}\|^2,
\]
whence $\|w_{k+1}\|=1$.

\smallskip

The same procedure works in the remaining subcases as well; for $k\in S_2$
the eigenvalues $\lambda_j$ and $\nu_j$ are simple, while for $k\in S_3$,
account must be taken of the multiplicity of $\vL$.

\end{proof}

\begin{remark}\label{torusremark}
The proof shows that, if $w_{k+1}$ is not required to be real, each
term $\alpha_ju_j$ is determined only up to a multiple
$\exp(\ei\theta_{k+1,j})$. Thus, the
possible polygons $\be$ corresponding to a given pattern $\gamma$
lie on a torus. The angle coordinates are studied in the next section.
\end{remark}

We conclude by making  a choice of functionally independent action variables.

\begin{definition}\label{defineindex}
Let ${\mathcal I}$ be the set of pairs $(i,j)$ satisfying
$1\leq i\leq n-1, 1\leq j\leq i$
which index eigenvalues $\laij$ such that $\laij$ is not forced to be $0$
or $\vL$,  with the further property that $\lambda_{i,j+1}$ is
not forced to be
$0$ (this last condition says that in each row we throw away the right-most
$j$ such that  $\laij$ is not forced to be $0$).
\end{definition}

\begin{corollary}
The set ${\mathcal I}$ indexes a functionally
independent set of action variables $\laij$.
\end{corollary}
\begin{proof}
Indeed, these action variables map onto a
polyhedron of dimension equal to the cardinality of ${\mathcal I}$.
\end{proof}

\begin{remark}
For general coadjoint orbits, one can define a complete set of constants
of motion that reduce to the $\alpha_j$ in the rank one case; the construction
also makes use of Gel'fand-Tsetlin patterns. Action variables which
generate $2\pi$-periodic flows are not known, however.
\end{remark}

\section{Angle Variables and Four-Point Functions}\label{angles}
In this section, we construct angle variables $\thij$ conjugate
to the action variables $\laij$ discussed thus far. The angles
are implicit in Corollary \ref{twopi} and Remark \ref{torusremark};
what we now find is a global description.

\subsection{Four-point functions and polygons}

The geometric picture in \cite{KapovichMillson2} serves as model.
For the moment, think of the sides $e_j$ as vectors in $\bbR^3$.
 The action
variables are the lengths of the diagonals $A_i=e_1+\ldots+e_{i+1}$
of the polygon. The corresponding conjugate
angle is the {\em oriented} dihedral angle between the two triangles spanned,
respectively, by $A_{i-1},e_{i+1},A_i$ and $A_i,e_{i+2},A_{i+1}$.
By this we
mean the oriented angle between the two normal vectors to the triangles.
These two vectors are elements of the plane orthogonal to $A_i$. We
orient this plane so that a positively oriented basis for the
plane followed by $A_i$ is a positively oriented basis for $\bbR^3$.

\begin{remark}
In an oriented plane $\Pi$ equipped with a positive definite inner product
$U\cdot V,$  we can define the oriented angle $\angle (U,V)$ for a pair of
vectors $U$ and $V$ in $\Pi$ as follows. First we say that two unit vectors
$U,V$ make an angle of ninety degrees if $U\cdot V=0$ and the basis $\{U,V\}$
is
positively oriented. We let $J$ be the operation of rotation by ninety degrees.
We make $\Pi$ into a complex vector space by defining $\imath V := JV$.
Then the unit circle in $\bbC$ acts simply-transitively on the oriented
lines in $\Pi$. We define $\angle (U,V)=\theta$ if
$\exp(\imath \theta) U$ is a positive real multiple of $V$.
If $\theta = \angle (U,V)$ then we have
\begin{eqnarray*}
\cos \theta &=& \frac{U\cdot V}{\|U\| \|V\|}\\
\sin \theta &=&  \frac{JU\cdot V}{\|U\| \|V\|}
\end{eqnarray*}
\end{remark}

For the case at hand, the  oriented angle $\theta_i$ is given by
\begin{align}
\cos\theta_i &=\ \  \frac{(A_i\times e_{i+1})\centerdot (A_i\times e_{i+2})}
    {\|A_i\times e_{i+1}\|\, \|A_i\times e_{i+2}\|} \label{cos}\\
 \sin\theta_i &=  \frac{(A_i\times e_{i+1})\times (A_i\times e_{i+2})
\centerdot A_i}
  {\|A_i\times e_{i+1}\|\, \|  A_i\times e_{i+2}\|\, \|A_i\|}.\label{sin}
\end{align}

Note that $\theta_i=0$ when the triangles are coplanar, so that
the collection of planar polygons forms a reference cross-section
for the angle variables.

We now transfer (\ref{cos}) and (\ref{sin}) back to our Lie algebra $\calH_2^0$
of tracefree Hermitean $2\times 2$ matrices.
Define $f:\bbR^3\to\calH_2^0$ by
\begin{equation}\label{translate}
f:\bx=(x_1,x_2,x_3)\mapsto \hx=\frac{1}{2}
  \begin{pmatrix}
   x_1 & x_2 +\ei x_3 \\
   x_2-\ei x_3 & -x_1
   \end{pmatrix}.
\end{equation}
Then $\widehat{\bx\times\by}=\ei[\hx,\hy]$,
$\bx\cdot\by=2\Tr\hx\hy$,
and a vector in the $x_3=0$ plane corresponds to a real symmetric
matrix. (Thus, a planar polygon is represented by a symmetric matrix,
cf. Theorem \ref{inversetheorem}).

We return to identifying vectors with matrices via (\ref{translate}).\\
Let $\lambda>0$ and $-\lambda$ be the eigenvalues of $A_i$, with
orthonormal eigenvectors $u,v$,
so that $A_i=\lambda(u\otimes u^*-v\otimes v^*)$.
Write, for notational simplicity,
$$e_{i+1}=r_1\wtens{1}-(r_1/2)\id, \quad
e_{i+2}=r_2\wtens{2}-(r_2/2)\id.$$
Then the numerator of (\ref{cos}) becomes (since
$\id$ does not contribute)
\begin{equation}\label{bigmess}
2\Tr\bigl(\ei[A_i,r_1\wtens{1}]
    \ei[A_i,r_2\wtens{2}]\bigr),
\end{equation}
and the numerator of (\ref{sin}) becomes
\begin{equation}\label{bigmess2}
2\|A_i\| \Tr \bigl(\ei A_i [r_1 \wtens{1},r_2 \wtens{2}]\bigr).
\end{equation}

\begin{definition}[\cite{BerceanuSchlichenmaier}]\label{fourpoint}
Let $a,b,c,d\in\bbCm$.
Define the {\em four-point function}  by
\[
F_4(a,b,c,d)=\frac{(a,b)(b,c)(c,d)(d,a)}{\|a\|^2\|b\|^2\|c\|^2\|d\|^2}
\]
where $(\cdot,\cdot)$ is the usual Hermitean inner product.
\end{definition}
Two properties of $F_4$ are important:
\begin{enumerate}

\item $F_4(a,b,c,d)$ may be thought of as function on $(\bbCPm)^4$; in
particular, $F_4$ is independent of the phases of its arguments.

\item $\overline{F_4(a,b,c,d)}=F_4(a,d,c,b)$ (plus other such symmetries).
\end{enumerate}

A longish calculation, using property (2),  gives the following.

\begin{proposition}\label{ReF_4} Expression (\ref{bigmess}) reduces to
\[16\lambda ^2 r_1r_2 \Real F_4(w_1,u,w_2,v).\]
Expression (\ref{bigmess2}) reduces to \[16\lambda ^2  r_1r_2 \Imag
F_4(w_1,u,w_2,v).\]
The denominator in (\ref{cos}) and (\ref{sin})
becomes
\[16\lambda ^2 r_1r_2 |F_4(w_1,u,w_2,v)|.\]
Thus, the oriented dihedral angle is $\theta=\arg F_4(w_1,u,w_2,v)$.
\end{proposition}

This formula, suitably adapted, will be shown to define the
conjugate angles in the more general case as well.

\smallskip

We mention, as an aside, that the argument of the four-point function
has an interesting geometric description.

\begin{theorem}\label{fourpointfunction}
Let $a_j,j=1,\ldots,4$ be four points in $\bbCm$ defining
points $p_j\in\bbCPm$. Construct a geodesic
quadrilateral $\pi$ in
$\bbCPm$ with vertices at the $p_j$. Let $\sigma$ be a two-chain with
boundary $\pi$ and
let $\omega$ be the K\"ahler form on $\bbCPm$.Then
\begin{equation}\label{integral}
\arg F_4(a_1,a_2,a_3,a_4)= - \int_{\sigma}\omega.
\end{equation}
\begin{proof}
Draw a geodesic segment (a diagonal of the quadrilateral) from $p_1$ to
$p_3$. The analogue of (\ref{integral})
for triangles was proved in \cite{HanganMasala}, see also
\cite[Ch. 7]{Goldman}.
Now choose $\sigma$ to be the union of
two two-chains each of which has as boundary one of the two
triangles created by drawing the diagonal $p_1p_3$. Combining
(\ref{integral})
for the triangles gives the  equation for the quadrilateral.
\end{proof}

\end{theorem}

\subsection{Construction of angle variables}

We will define the angle variables as in Proposition \ref{ReF_4}, via the
four-point function of the $w$'s associated with two consecutive edges
and eigenvectors of the diagonal between them. These vectors all involve
a choice of phase, and the first goal will be to remove the ambiguity.

Let $\Mr^0$ be the open subset of $\Mr$ on which
the interlacing inequalities $\laij>\lambda_{i-1,j}>\lambda_{i,j+1}$
are strict, and let $\Mrtw^0$ be its inverse image in $\Mrtw$.
We consider only polygons in $\Mrtw^0$, so that the (unconstrained)
eigenvalues and eigenvectors may be taken to be locally smooth
functions of $\be$.

Let $\phi^t$ be one of the $\laik$-flows defined in
Proposition \ref{phi-ij}.
We will follow the transformed polygon $\phi^t(\be)$.
Its $\ell$-th edge, $r_\ell w_\ell^t\otimes (w^t_\ell)^*$,
and the normalized $\laij$-eigenvector, $u_{ij}^t$, of the
diagonal $\phi^t(A_i)$, will depend on time $t$.
They may be taken to be
locally smooth on $\Mr^0$, but will depend on an initial choice,
while the polygon $\phi^t(\be)$ itself is well defined.

\begin{definition}\label{anglevariables}

Make smooth local choices of $w_\ell$ and $u_{ij}$. Here $u_{ij}$
is a (choice of) unit length eigenvector belonging to the eigenvalue
$\lambda_j$ of $A_i$. Define
$$\alpha_{ij}:\Mrtw^0\to\bbC,(i,j)\in \mathcal{I}, {\mbox {\ by\ }}
 \alpha_{ij}:\be\mapsto (w_{i+1}(\be),u_{ij}(\be))
(u_{ij}(\be),w_{i+2}(\be));$$
this depends on the phases of $w_{i+1},w_{i+2}$. (We will usually drop
the argument $\be$). Set
$$\beta_{ij}=  F_4(w_{i+1},u_{ij},w_{i+2},u_{i,j+1})=\alpha_{ij}
\overline{\alpha_{i,j+1}}.$$
The $\beta_{ij}$ are {\em independent} of all phase choices.
Finally, we define the angle variables $\thij,(i,j)\in \mathcal I $,  by
$$\thij=\arg \beta_{ij}.$$
\end{definition}

Clearly the number of four-point functions $\beta_{ij}$ is the same
as the number of independent, unconstrained $\laij$'s,
since for every $i$ there is one more $\laij$ than $\beta_{ij}$
and there are no $\beta_{ij}$'s corresponding
to the eigenvalues $0$ and $\vL$. Thus we obtain the correct
formal count of angle variables. We now prove
that the angle variables are well-defined on $\Mr^0$.

\begin{lemma}\label{nonvanishing}\hfill
\begin{enumerate}

\item All $|\alpha_{ij}|^2$ are constant under all
bending flows $\phi_{k\ell}$.
\item All $|\alpha_{ij}|^2$ are nonzero on $\Mrtw^0$.
\end{enumerate}
In particular, $\arg \beta_{ij}=\arg\alpha_{ij}\overline{\alpha_{i,j+1}}$
is defined.
\end{lemma}
\begin{proof}The first statement follows from Proposition \ref{detratio} and
Corollary \ref{rationalfunctions}. Indeed,
$$A_{i-1}=A_i - r_{i+1}\wtens{i+1}.$$
Hence $|(w_{i+1},u_{ij})|^2$, being a rational
function of action variables, is a constant of motion.
Likewise,
$$A_{i+1}=A_i+r_{i+2}\wtens{i+2}$$
implies that $|(w_{i+2},u_{ij})|^2$ is a constant of motion.

To prove the second statement we apply the Weinstein-Aronszajn formula to
obtain

\begin{equation}
  \frac{\det(z\id-A_{i-1})}{\det(z\id-A_i)}=1+r_{i+1}\sum_{j=1}^{m+1}
  \frac{|(w_{i+1},u_{ij})|^2}{z-\laij}.
\end{equation}

If $|\alpha_{ij}|=0$ then either $|(w_{i+1},u_{ij})|= 0$ or
$|(w_{i+2},u_{ij})|= 0$. Assume first that $|(w_{i+1},u_{ij})|= 0$.
>From the Weinstein-Aronszajn formula we see that if follows that
$\laij$ is not a
pole, so the $(z-\laij)$ in the
denominator of the left-hand side must cancel with one of the terms
in the numerator.
Hence one of the interlacing inequalities between the $i^{\text{th}}$ and
$(i-1)^{\text{st}}$ rows is not
strict, contradicting the assumption that $\be\in \Mrtw^0$. Similarly,
$(w_{i+2},u_{ij})\ne0$.
\end{proof}

\begin{lemma}\label{betaij-invariant}
The $\beta_{ij}$ are invariant under conjugation.
\end{lemma}
\begin{proof}
Let $g\in\Um$ and consider the conjugated polygon $g\be g^{-1}$. Its
$\ell^{\text{th}}$ edge is $r_\ell (gw_\ell)\otimes (gw_\ell)^*$.
However, the choice
$w_\ell(g\be g^{-1})$ made in Definition \ref{anglevariables} may not
coincide with $gw_\ell$. If they differ, it is by a multiple of modulus
one. The four-point function $\beta_{ij}$ is
not affected by such a factor.
\end{proof}

In the following we will make essential use of

\begin{remark}\label{geg-inv}
In view of the proof of Lemma \ref{betaij-invariant}, we may  replace
$w_\ell(g\be g^{-1})$
by $g w_\ell$ in calculations involving $\beta_{ij}$, and for the same
reason, $u_{ij}(g\be g^{-1})$ by $gu_{ij}$.
\end{remark}

We will now compute the Poisson brackets of the
action variables with the angle variables.

\begin{lemma}\label{poissonbracketsofactionandangle}

$$\{\lail,\thij\}=
\begin{cases}
  &\phantom{-}1,\, l= j\\
  &-1,\,l=j+1\\
  &\phantom{-}0,\, l\ne j,j+1
\end{cases}$$
\end{lemma}

\begin{proof}
We will verify, using (\ref{phi-ij}), that
$$\beta_{ij}(\phi_{il}^t(\be))=
\begin{cases}
\phantom{\exp(-\ei t)\,}\beta_{ij}(\be),\,&l\ne j,j+1,\\
\phantom{-}\exp(\ei t)\,\beta_{ij}(\be),\,&l=j,\\
\exp(-\ei t)\,\beta_{ij}(\be),\,&l=j+1.
\end{cases}$$
Note from (\ref{phi-ij}) that the $i^{\text{th}}$ diagonal $A_i$
of $\be$ and the ($i+2$)-nd edge
are fixed under $\phi_{il}^t$. Hence the normalized eigenvectors
$u_{ij}$ of $A_i$ are also fixed.
Now abbreviate $g_t=\exp(\ei tE_l(A_i))$, and as explained in
Remark \ref{betaij-invariant},
make the replacement
$$w_{i+1}(\phi_{il}^t(\be))= w_{i+1}(g_t \be g_t^{-1})\rightsquigarrow
g_tw_{i+1}(\be).$$

We obtain
 \begin{align*}
\beta_{ij}(\phi_{il}^t(\be))&= (g_tw_{i+1},u_{ij})(u_{ij},w_{i+2})
(w_{i+2},u_{i,j+1})(u_{i,j+1},g_t w_{i+1})\\
&= (w_{i+1},g_t^{-1} u_{ij})(u_{ij},w_{i+2})
(w_{i+2},u_{i,j+1})(g_t^{-1} u_{i,j+1},w_{i+1}).
\end{align*}

Since $E_l(A_i)u_{ij}= \delta_{jl}u_{ij}$ the lemma follows by definition
of $g_t$.
\end{proof}

\begin{lemma}

$$\{\laij,\theta_{kl}\}= 0,i\ne k.$$

\end{lemma}

\begin{proof}

If $i<k$ then  the $k^{\text{th}}$ diagonal, the ($k+1$)-st edge, and the ($k+2$)-nd
edge are fixed by
the bending flow $\phiij^t$, and hence $\theta_{kl}$ is unchanged.

If $i>k$, then the $k^{\text{th}}$ diagonal,the ($k+1$)-st edge and the ($k+2$)-nd
edge are rigidly moved by
the $g_t$ under the bending flow $\phiij^t$, and
hence $\theta_{kl}$ is unchanged. (Note that
Remark \ref{betaij-invariant} is used once more).
\end{proof}

To remove the redundancy in the $\laij$, we define new action variables
$\mu_{ij}$ by the formula
\begin{equation}\label{newactions}
\mu_{ij}= \sum_{k=1}^j \lambda_{ik}.
\end{equation}

\noindent As a consequence of the two preceding lemmas we obtain

\begin{proposition}
The action variables $\{\mu_{ij}\}$ and the angle variables $\{\thij\}$
are conjugate
$$\{\mu_{ij},\theta_{kl}\}=
\begin{cases}
&1,\,i=k,j=l\\
&0,\,otherwise.
\end{cases}$$

\end{proposition}

We deduce two corollaries.

\begin{corollary}
The angle variables are functionally independent.
\end{corollary}

\begin{corollary}\label{movelevels}
The Hamiltonian flows of the new action variables $\{\mu_{ij}\}$ permute
the simultaneous level sets $\{\thij = c_{ij}, (i,j)\in\mathcal{I}\}$
transitively.
\end{corollary}

\medskip
We now begin the proof that
$$\{\thij,\theta_{kl}\}=0.$$

Recall that $\calSm$ is the space of real symmetric
$(m+1)\times(m+1)$ matrices. Let $\sigma:\calHm\to \calHm$ be complex
conjugation. Then $\calSm$ is the fixed subspace of $\sigma$. The
following lemma is immediate from (\ref{bracket}):

\begin{lemma}
The involution $\sigma$ is anti-Poisson (a Poisson isomorphism from
$\calHm$ equipped with the Lie Poisson tensor to $\calHm$ equipped
with the negative of the Lie Poisson tensor).
\end{lemma}

We obtain

\begin{corollary}
If $f$ and $g$ are constant on $\mathcal S_{m+1}$, then $\{f,g\}$ vanishes on
$\mathcal S_{m+1}$.
\end{corollary}
\begin{proof}
Let $\pi(.,.)$ be the Lie Poisson bivector considered as a
skew-symmetric bilinear form on the
cotangent bundle of $\mathcal H_{m+1}$.
For $x\in \mathcal S_{m+1}$
and $u,v$ cotangent vectors at $x$, the Lemma gives
$\pi_x(u,v)=-\pi_x(\sigma u,\sigma v)$.
If $u$ and $v$ are conormal covectors at $x$ then they are in
the $(-1)$-eigenspace for $\sigma$,
and therefore $\pi_x(u,v)=0$. But if $f$ and $g$ are constant on
$\calSm$, then $df_x$ and $dg_x$ are conormal at $x$.
\end{proof}

As an immediate consequence we have
\begin{lemma}\label{vanishingofpoissonbrackets}
If $f$ and $g$ are constant on $M_{\br}(\calSm)$, then $\{f,g\}$ vanishes on
$M_{\br}(\calSm)$.
\end{lemma}

Our next goal is to prove that the simultaneous zero level set of
the angle variables is
$M_{\br}(\calSm)$. In order to obtain this we will need
two technical lemmas to
handle the regions $S_1$ and $S_3$ (in the notation of
Proposition \ref{bigPdim}).
The first lemma will be used to deal with the region $S_3$.

\begin{lemma}\label{region3}
Let $V_i=\ker(A_i-\vL\id), n-m-2\le i\le n-1$. Then
$$V_{n-1}\supset V_{n-2}\supset \cdots \supset V_{n-m-2}=\{0\}.$$
Moreover (recalling $A_i=A_{i-1}+ r_{i+1} \wtens{i+1}$) we have
$$V_{i-1}=\{v\in V_i: (v,w_{i+1})=0\}.$$
\end{lemma}

\begin{proof}
Let $v\in V_{i-1}$ and $\|v\|=1$. Then
$$\vL = (A_{i-1}v,v)= (A_iv,v)-r_{i+1}|(w_{i+1},v)|^2.$$
But $\vL$ is the largest eigenvalue of $A_i$ so $(A_iv,v)\le \vL$.
Hence the above equation can hold if and only if
$$(A_iv,v)=\vL \ \  (\mbox{so}\,  v\in V_i) \ \ \mbox{and}\ \ (w_{i+1},v)=0.$$
\end{proof}

\begin{corollary}
Let $w_{i+1}^{\vL}$ be the orthogonal projection of $w_{i+1}$ on
the $\vL$-eigenspace of $A_{i-1}$. Then
$$w_{i+1}^{\vL}=0.$$
\end{corollary}
The next lemma will be used to deal with the region $S_1$.

\begin{lemma}\label{region1}
Let $U_i=\ker A_i, 1\le i\le m$. Then
$$U_1\supset U_2\supset \cdots \supset U_m=\{0\}.$$
Moreover
$$U_i=\{u\in U_{i-1}: (u,w_{i+1})=0\}.$$

\end{lemma}

\begin{proof}
Suppose $A_iu=0$. Then
$$0=(A_iu,u)=(A_{i-1}u,u)+r_{i+1}|(w_{i+1},u)|^2.$$
But $A_{i-1}$ is positive semidefinite and $r_{i+1}>0$.
Hence $u\in \ker A_{i-1}$ and $(u,w_{i+1})=0$.
\end{proof}

\begin{corollary}
Let $w_{i+1}^0$ be the projection of $w_{i+1}$ on $\ker A_i$. Then
$$w_{i+1}^0=0.$$
\end{corollary}

Now we can prove the result we need. Let $Z(\Theta)$ be the
simultaneous zero level set of the angle variables $\{\thij\}$.

\begin{proposition}\label{zerolevel}
$$Z(\Theta)= M_{\br}(\mathcal S_{m+1}).$$
\end{proposition}

\begin{proof}

The inclusion
$$ M_{\br}(\calSm)\subset Z(\Theta)$$
is obvious (all the edges and diagonals are real, so the
eigenvectors are real, so the
$\beta_{ij}$ are real). The point is to prove the reverse inclusion.
We will assume $n\ge 2(m+1)$ and leave the case $n\le 2m+1$, which is
similar, to the reader.

Given a polygon $\be$ with all $\thij=0$. We wish to show that
a sequence of conjugations of $\be$ by elements of $\Um$ will make
all sides $e_k$ real symmetric, or equivalently, all the $w_k$ real.
The proof is by descending induction, starting with
the last diagonal $A_{n-1}=e_1+\cdots+e_n=\vL \id$, which is
of course real symmetric. First, conjugate $\be$ by $g\in \Um$
(without changing $A_{n-1}$) to arrange that $A_{n-2}$ is diagonal,
hence real. This moves all the $w_k$ to $gw_k$, but in the sequel
we do not need to keep track of those changes.
Now we know that $A_{n-3}$ has the form
$$A_{n-3}=A_{n-2}-r_{n-1}\wtens{n-1},$$
and we want to show that we can move $w_{n-1}$ to a real vector. We have
$$\ker (A_{n-2}-\vL\id) = \{\epsilon_1,\ldots,\epsilon_m\},$$
where $\{\epsilon_1,\ldots,\epsilon_{m+1}\}$ is the standard basis for $\bbCm$.
Suppose $A_{n-2}\epsilon_{m+1}= \mu \epsilon_{m+1}, \mu=\vL-r_n$.

Write $w_{n-1}$ in the form
$w_{n-1}=w_{n-1}^\vL+w_{n-1}^\perp$, where $w_{n-1}^\vL$ is the orthogonal
projection of $w_{n-1}$ onto $\ker(A_{n-2}-\vL\id)$. Hence there
exists $z\in \bbC$ such
that  $w_{n-1}^\perp = z \epsilon_{m+1}$. Since $w_{n-1}$ is defined
only up to a complex
multiple of unit length, we may multiply $w_{n-1}$ by an
element of $S^1$ in order to
arrange that $z$ be real. Let $c = \|w_{n-1}^\vL\|$.
Now choose $g\in \Um$ such that
$g\epsilon_{m+1}=\epsilon_{m+1}$ and $gw_{n-1}^\vL = c \epsilon_m$.Then
$gA_{n-2}g^{-1}= A_{n-2}$ (because $g\epsilon_{m+1}=\epsilon_{m+1}$ and
$gw_{n-1}= c\epsilon_m+z\epsilon_{m+1})$. We change $\be=(e_1,\ldots,e_n)$ to
$g\be g^{-1}=(ge_1g^{-1},\ldots,ge_ng^{-1})$.

Next, we show how to find a conjugation $g\be g^{-1}$ that keeps
$A_{n-2},A_{n-3}$ and $w_{n-1}$ real and also makes $gw_{n-2}$ real.
This step exhibits the general pattern.

By Lemma \ref{region3},
\begin{align*}
\ker (A_{n-3}-\vL\id)&=\{v\in \ker(A_{n-2}-\vL\id):(v,w_{n-1})=0\}\\
&= \spann \{\epsilon_1,\ldots,\epsilon_{m-1}\}.
\end{align*}
The matrix $A_{n-3}$ has two new eigenvalues (in addition to $\vL$);
let their eigenvectors be $u_{n-3,m+1},u_{n-3,m}$.
There is one angle variable
\begin{align*}
\theta_{n-3,m+1}=&\arg [(w_{n-2},u_{n-3,m})(u_{n-3,m},w_{n-1})\\
&(w_{n-1},u_{n-3,m+1})(u_{n-3,m+1},w_{n-2})]
\end{align*}
We have seen that $A_{n-3}$ is real symmetric, hence $u_{n-3,j}$
can be chosen to be real for all
$1\le j \le m+1$. Since $w_{n-1}$ is real, we may normalize
$u_{n-3,m}$ and $u_{n-3,m+1}$ so
that $(w_{n-1},u_{n-3,m})>0$ and $(w_{n-1},u_{n-3,m+1})>0$.
Since, by assumption,
$\theta_{n-3,m+1}=0$, we have
$$\arg (w_{n-2},u_{n-3,m+1})=\arg (w_{n-2},u_{n-3,m}).$$
Hence by multiplying $w_{n-2}$ by an element in $S^1$ we may
assume that $(w_{n-2},u_{n-3,m+1})$
and $(w_{n-2},u_{n-3,m})$ are real. Now we may write
$$w_{n-2}=w_{n-2}^{\vL}+w_{n-2}^{\perp},$$
where
$$w_{n-2}^{\vL}\in ker((A_{n-3}-\vL\id)=
\spann\{\epsilon_1,\ldots,\epsilon_{m-1}\}$$
and
$$w_{n-2}^{\perp}\in \spann \{u_{n-3,m},u_{n-3,m+1}\}
= \spann\{\epsilon_m,\epsilon_{m+1}\}.$$
We have arranged for $w_{n-2}^{\perp}$ to be real.
Choose $g\in \Um$ with
$g\epsilon_m=\epsilon_m$ and $g\epsilon_{m+1}=\epsilon_{m+1}$ such that
$$g w_{n-2}^\vL= c'\epsilon_{m-1},$$
with $c'=\|w_{n-2}^\vL\|$ as in the preceding step.
Now change $\be$ to
$g\be g^{-1}$ and proceed to $w_{n-3}$.

We continue in this way until $\ker(A_k-\vL\id)=0$ and we
enter the region $S_2$. The argument for this region is simpler and
is left to the reader. Note that the
vanishing of the angle variables says that {\em all}
the coordinates $(w_k,u_{k-1,j})$ in the eigenvector
basis of $A_{k-1}$ have a common phase which can be
eliminated by multiplication by
an element of $S^1$; no conjugation is needed, so the preceding
edges all remain real symmetric. However, the zero eigenvalue, which is
unavoidable when we enter region $S_1$, causes new problems,
and Lemma \ref{region1} is required.

Suppose then we have proved that $A_m$ is real (note that
$A_m$ has rank $m$). We want to
prove that $A_{m-1}$ is real. We know that
$$A_m=A_{m-1}+ r_{m+1}\wtens{m+1},$$
and since $\ker A_m=\{0\}$,
we have enough angle variables to prove that all
coordinates of $w_{m+1}$ have a common
phase. We clear this phase as before and move on to $A_{m-2}$.
We have
$A_{m-1}=A_{m-2}+r_m \wtens{m}$,
and wish to prove that one can make
$w_m$ real without destroying reality of $w_{n-1},\ldots,w_{m+1}$.
Write $w_m = w_m^{\perp}+w_m^0$ with $A_{m-1}w_m^0=0$ and $w_m^{\perp}$
orthogonal to
$\ker A_{m-2}$ (the latter has dimension $2$). By the
corollary to Lemma \ref{region1},
we have $w_m^0=0$. Also, we have enough
angle variables to conclude that the coordinates of $w_m^{\perp}$
relative to the eigenvectors
of $A_{m-1}$ orthogonal to $\ker A_{m-2}$ have a common phase.
Thus, no conjugations are required to make $w_m$ real,
and all preceding edges remain real symmetric. Now continue.
\end{proof}

We remark that the proof could equally well be done by ascending
induction; in that case, region $S_1$ would be the one requiring
conjugations, while an overall scaling would do in $S_2,S_3$.

\medskip
We are now ready to prove
\begin{proposition}
$$\{\thij,\theta_{kl}\}=0.$$
\end{proposition}

\begin{proof}
Let $\be\in \Mr$ be given. By Corollary \ref{movelevels},
the bending deformations flows permute the
level sets of the $\theta_{ij}$'s transitively. Hence we may
apply a bending $\phi$ to move
$\be$ into $Z(\Theta)$. Since $\phi$ is symplectic and the
Hamiltonian vector fields of the $\thij$
are invariant under bending, we have
$$\{\thij,\theta_{kl}\}(\be)= \{\thij,\theta_{kl}\}(\phi\be).$$
But by Proposition \ref{zerolevel}
$$Z(\Theta)= M_{\br}(\mathcal S_{m+1}).$$
Hence by Lemma \ref{vanishingofpoissonbrackets}
$$\{\thij,\theta_{kl}\}=0.$$
\end{proof}

\section{The duality between the bending systems and the Gel'fand-Tsetlin
systems
on Grassmannians}

In this section we use Gel'fand-MacPherson duality, following
\cite{HausmannKnutson}
for the case of $m=1$, to
show that the bending system is equivalent to the Gel'fand-Tsetlin integrable
system
(as defined in \cite{GuilleminSternberg}) on a torus quotient of the
Grassmannian
$G(m+1,\bbC^n)$. This equivalence will explain the  appearance and form of
the Gel'fand-Tsetlin patterns in \S 8.

Our first goal is to construct a symplectomorphism $\Phi$ from $\Mr$ to
a symplectic quotient of $G(m+1,\bbC^n)$ by the
$n$-torus $T$ of diagonal matrices in $\Un$.
This is the symplectic version of Gel'fand-MacPherson duality.

Let $\cm$ denote the vector space of $(m+1)\times n$ complex
matrices. We give $\cm$ the Hermitean form $(\ ,\ )$ defined by
$(X,Y)=2\Tr(X^*Y)$, and thus $\cm$ is a symplectic vector space. The  product
group $\Um\times\Un$ acts isometrically and symplectically on $\cm$.  Denote
the $i^{\text{th}}$ row
(resp. $j^{\text{th}}$ column) of $N\in \cm$  by $R_i$ (resp.
$C_j$).
\begin{pproposition} The action of $\Un$ has momentum map
$$\mu_{\Un}:\cm\to\calH_n,\quad \mu_\Un:N\mapsto N^*N.$$
In particular, the momentum map for the $T$-action is
$$\mu_T:N\mapsto(\|C_1\|^2,\ldots,\|C_n\|^2).$$
The momentum map for the $\Um$ action is
$$\mu_\Um:\cm\to\calHm,\quad \mu_\Um:N\mapsto NN^*.$$
\end{pproposition}

Note that
\[
\mu_\Um(N)=\sum_{i-1}^n C_i\otimes C_i^*.
\]
This will provide the connection with polygons.

We construct the desired symplectomorphism by computing the
symplectic quotient corresponding to the $\mu_T$-level $\br$ and the
$\mu_\Um$ level
$\vL \id$ in two different orders. If we first quotient with respect to $T$
with momentum level $\br$ and then with respect to $\Um$ with
momentum level $\vL\id$,
we get the space $\Mr$. In order to see this, we note that the (left) action
of $\prod_1^n \Um$ on $\cm$ (acting on the columns) commutes with the
(right) action of
$T$ (in fact one obtains a dual pair in the sense of Howe, see
\cite{KazhdanKostantSternberg}). We first compute the symplectic
quotient by $T$.

\begin{llemma}
\begin{enumerate}
\item The momentum map $\mu_{(\Um)^n}$ induces an embedding of the
symplectic quotient
$\mu_{T}^{-1}(\br)/T$ into $\Pi_1^n\calH_{m+1}$, with image
$\prod_1^n\orb_{r_i}$.
\item The form on $\mu_{T}^{-1}(\br)/T$ induced by reducing the form $2 \Imag\,X^*
Y$
when carried over to $\prod_1^n\orb_{r_i}$ agrees with the Kostant-Kirillov
form $\omega_{KK}$.
\end{enumerate}
\end{llemma}

\begin{proof} The first statement follows because it is a general feature of dual pairs, see
\cite{KazhdanKostantSternberg}, that the
momentum map for one action embeds the symplectic quotient of the other
as an orbit in
(the dual of) the Lie algebra of the first group. This principle, applied to the
pair
$(\Um)^n\times T$, implies the first statement in the lemma. The second follows
from a straight-forward computation.
\end{proof}
Thus we have identified the quotient by $T$ with the correct product of rank
one orbits
in $\calHm$. Clearly, after taking the symplectic quotient of this product
by the diagonal action of $\Um$ (at momentum level $\vL\id$), we obtain $\Mr$.

Suppose instead we first quotient with
respect to $\Um$ and momentum level $\vL\id$. We get the Grassmannian
$G(m+1,\bbC^n)$ with a certain $\Un$-invariant symplectic structure.

\begin{llemma}
The momentum map $\mu_{\Un}$ induces an embedding of the symplectic quotient
$\mu_{\Um}^{-1}(\vL \id)/\Um$ into $\calH_n$, with image
the $\Un$-orbit $\orb_{\vL}$ consisting of those matrices that
have eigenvalue $\vL$ with multiplicity $m+1$ and eigenvalue $0$ with
multiplicity
$n-m-1$.
\end{llemma}
\begin{proof}The argument is the analogous to the previous case, only
this time we use the dual pair $\Um\times \Un$.
\end{proof}

Denote the torus quotient at momentum level $\br$ of the Grassmannian
with the Kostant-Kirillov symplectic structure corresponding to $\vL$
by $\cm_{\vL}$. We have now obtained the desired symplectomorphism $\Phi$
from $\Mr$ to $\cm_{\vL}$.

Of course this symplectomorphism gives a Poisson isomorphism between the
Poisson algebras of smooth functions of $\Mr$ and $\cm_{\vL}$.
However, we want to make this more explicit and to localize it.
Let $\cm_{\br,\vL}$ be
the subset of $\cm$ consisting of matrices $N$ such that $\|C_j\|^2=r_j$ and
$N^*N=\vL \id$. Thus we have $\Um\times T$ quotient mappings
$\pi_1:\cm_{\br,\vL}\to \Mr$
(first quotient by $T$ then by $\Um$) and $\pi_2:\cm_{\br,\vL}\to \cm_{\vL}$
(first quotient by $\Um$ then by $T$). We use the
mappings $\pi_1$ and $\pi_2$ to realize (and localize) the Poisson
isomorphism $\Phi$ from above.
Let $f$ be a function which is smooth on an open subset of $\Mr$. Use $\pi_1$ to
pull $f$ back to a $\Um\times T$-saturated open subset of $\cm_{\br,\vL}$.
Since $\pi_2$ is a quotient map and $\pi_1^*f$ is invariant under $\Um$,
we can first descend it to
to a $T$-saturated open subset of the Grassmannian, then to the torus quotient
of
that open set, which is an open subset of $\cm_{\vL}$. We note that $\Phi$ is
determined by the equation
$$\Phi(\pi_1(N))= \pi_2(N).$$

We now briefly review the Gel'fand-Tsetlin integrable system - for the details
see \cite{GuilleminSternberg}. As before, we identify the space $\calH_n$
of $n\times n$ Hermitean matrices with the dual of the Lie algebra of $\Un$.
We now
construct $n(n+1)/2$ Poisson commuting functions on $\calH_n$ which are smooth
on a dense open subset. Let $X\in \calH_n$. Let $\beta_i(X)$
be the principal
$i\times i$ diagonal block. Define $\gamma_{ij}$ on $\calH_n$ by
$$\gamma_{ij}(X) = \lambda_j(\beta_i(X)),$$
where $\lambda_j$ is the $j^{\text{th}}$ eigenvalue of the block. As usual, we assume
that
the eigenvalues of the $i^{\text{th}}$ block are arranged in nonincreasing order.
It is proved in \cite{GuilleminSternberg}
that the $\gamma_{ij}$'s Poisson commute. We note that the $\gamma_{nj}$ are
Casimirs.
The restrictions of the remaining Gel'fand-Tsetlin Hamiltonians to generic
orbits
are functionally independent and give rise to  integrable system on such
orbits.
The eigenvalues of the blocks interlace and can be arranged in a
``Gel'fand-Tsetlin" pattern, shown here for  $n=6$.
\[
\begin{array}{ccccccccccc}
\gamma_{61}&{}&\gamma_{62}&{}&\gamma_{63}&{}&\gamma_{64}&{}&\gamma_{65}&{}&\gamma_{66}\\
{}&\gamma_{51}&{}&\gamma_{52}&{}&\gamma_{53}&{}&\gamma_{54}&{}&\gamma_{55}&{}\\
{}&{}&\gamma_{41}&{}&\gamma_{42}&{}&\gamma_{43}&{}&\gamma_{44}&{}&{}\\
{}&{}&{}&\gamma_{31}&{}&\gamma_{32}&{}&\gamma_{33}&{}&{}&{}\\
{}&{}&{}&{}&\gamma_{21}&{}&\gamma_{22}&{}&{}&{}&{}\\
{}&{}&{}&{}&{}&\gamma_{11}&{}&{}&{}&{}&{}
\end{array}
\]
\centerline{Figure 3}

Since we are dealing with a degenerate orbit here (the Grassmannian), many of
the
$\gamma_{ij}$'s (at the ends of the rows) will be zero (see
Remark \ref{zerogts}  below, and
Figure~2 above).
The next proposition, combined with the earlier sections, shows how to extract a
functionally independent set of Gel'fand-Tsetlin Hamiltonians and obtain
angle variables
for the Gel'fand-Tsetlin system on the Grassmannian.

\begin{pproposition}\label{GTsEqualsEivals}
$\Phi^*\gamma_{ij}=\lambda_{ij}$.
\end{pproposition}
\begin{proof}
Let $\id_k$ be the diagonal
$n$ by $n$ matrix whose first $k$ eigenvalues are equal to $1$ and last $n-k$
eigenvalues are equal to $0$. We use $\id_k$ to ``truncate" $N, N^*N$ and
$NN^*$.
Put $N_k:= N\id_k $. Then
\begin{align*}
\mu_{\Un}(N_k)= &\id_k N^*N \id_k\\
\mu_{\Um}(M_k)=& N\id_k \id_k N^*.
\end{align*}
The matrix on the first line is $\beta_k(N^*N)$, the
principal $k$ by $k$ block of the $n\times n$ matrix $N^*N$, and the matrix
on the second line is the diagonal
$A_{k-1}= C_1 C_1^* + C_2C_2^* +\cdots + C_kC_k^*$.
The matrices $\id_k N^*N \id_k$ and $N \id_k \id_k N^*$ have the same nonzero
eigenvalues. But the eigenvalues of the second matrix are the bending
Hamiltonians
$\lambda_{kj}$, and the eigenvalues of the first matrix are the
Gel'fand-Tsetlin
Hamiltonians $\gamma_{kj}$.  Finally we observe that

\begin{align*}
\gamma_{ij}(\Phi(\pi_1(M)))&=\gamma_{ij}(\pi_2 (M))=\lambda_j(\beta_i(\pi_M))\\
&=\lambda_j(A_i(\pi_1(M)))=\lambda_{ij}(\pi_1(M)).
\end{align*}
\end{proof}

We conclude this section with three remarks.

\begin{rremark}\label{zerogts}
Proposition \ref{GTsEqualsEivals}
explains the appearance of Gel'fand-Tsetlin patterns
in connection with the bending Hamiltonians. The appearance of the zeroes at
the
end of the rows in our patterns is explained because the Gel'fand-Tsetlin
system
in question is defined on a subset of the Hermitean matrices of rank at most
$m+1$.
Hence $\gamma_{ij}=0,j>m+1$.
\end{rremark}

\begin{rremark}
The reconstruction process in \S 9 may be interpreted as saying that
the class of patterns
introduced in \S 8 is precisely the class corresponding to Hermitean matrices
of the form $N^*N$, where $N$ is as above.
\end{rremark}

\begin{rremark}
Fixing the row sums in the patterns in \S 8 to be
partial sums of the $r_j$ corresponds to taking the
quotient of the Grassmannian by $T$ (at level $\br$).
\end{rremark}

\section{Pieri's formula and the duality at the quantum level}

In this section we will assume that the $r_i$'s are (positive) integers
and that $\Lambda = (r_1+\cdots + r_n)/(m+1)$ is an integer.The
orbit $\orb_{r_i}$ then corresponds under geometric quantization to
the irreducible representation $\mathcal{S}^{r_i}(V)$ of $\Um$, where
$V$ denotes the standard (or vector) representation of
$\Um$ on $\bbCm$ and $\mathcal{S}^{r_i}(V)$ the $r_i^{\text{th}}$ symmetric
power.

The (classical) duality result of the last section should have a
quantum version.
We note that the duality connected an integrable system (bending)
on a symplectic
quotient of $\prod_1^n \orb_{r_i}$ by the diagonal action of $\Um$
and an integrable
system (Gel'fand-Tsetlin) on a torus quotient of the
Grassmannian $G(m+1,n)$.
Thus, according to geometric quantization, at the quantum level we would expect
a relation between an $n$-fold tensor product multiplicity for $\GLm$
and a weight
multiplicity for a Cartan power of the the $(m+1)^{\text{st}}$ exterior
power of $\GLn$.
The bending system provides a (singular) real polarization of the space
$\Mr$, the symplectic quotient (at level $\vL\id$) of
$\prod_i \orb_{r_1}$. Thus the
number of lattice points in the momentum polyhedron
${\bf P}$ for bending should be
equal to the multiplicity of the the $1$-dimensional
representation $(det)^{\vL}$ in
$\otimes_1^n \mathcal{S}^{r_i}(V)$. But on the other hand,
the Gel'fand-Tsetlin system
is a real polarization of the torus quotient of the
Grassmannian (at level $\br$) where
the Grassmannian is given the symplectic structure which corresponds to the
orbit
of $\Un$ through the diagonal matrix with $m+1$ $\vL$'s
and $n-m-1$ zeroes. Thus the above number of lattice points should also
be the multiplicity of the $\br^{\text{th}}$ weight space in
$C^{\vL}\bigwedge^{m+1}V$,
the $\vL^{\text{th}}$ Cartan power of the $m+1$-st exterior power of the
vector representation
$V$ of $\GLn$. (We recall that if $W^\nu$ is a
representation with highest weight $\nu$,
then the $p$-th Cartan power
$C^pW^\nu$ is the irreducible representation with highest weight
$p\nu$). This equality of multiplicities predicted is in fact true, and will be
proved
below.

\begin{rremark}
It is unfortunate that the theory of geometric quantization
using a real polarization
is not sufficiently well developed to allow  one to deduce
theorems in representation
theory from equalities of numbers of lattice points in
momentum polyhedra. At this time
we can only regard such equalities as predictions of theorems in
representation theory.
\end{rremark}

We first note how the interlacing of the spectra of the perturbed
matrix and the unperturbed matrix (see \S 8)from the
Weinstein-Aronszajn formula
predicts Pieri's formula in representation theory.

\subsection{The Weinstein-Aronszajn and Pieri formulas}
We recall Pieri's formula for tensoring an  irreducible polynomial
representation of $\Um$ with a symmetric power of the vector
representation \cite[\S A.1]{FultonHarris}.

\begin{theorem}[Pieri's Formula]
Let $\lambda=(\lambda_1,\ldots,\lambda_{m+1})$ be the highest weight
of the  polynomial representation $V(\lambda_1,\ldots,\lambda_{m+1})$ of
$\Um$. Let $k$ be a positive integer. Then
$$ V(\lambda_1,\ldots,\lambda_{m+1})\otimes \mathcal{S}^{k}(V)=
\oplus V(\nu_1,\ldots,\nu_{m+1})$$
where the sum is taken over all dominant $\nu = (\nu_1,\ldots,\nu_n)$
satisfying
$$\nu_1\ge \lambda_1 \ge \nu_2 \ge \ldots \ge \nu_{m+1} \ge \lambda_{m+1} \ge
0$$
and
$$\sum_{i=1}^{m+1} \nu_i = \sum_{i=1}^{m+1} \nu_i + k.$$
\end{theorem}

This is of course Proposition \ref{interlace} restricted to integer
eigenvalues.
If $A\in \mathcal{O}_{\lambda}$, then the spectrum of the rank one
perturbation
$A+k\  w\otimes w^*$ satisfies the interlacing and row sum
conditions of the Pieri formula.

\subsection{Duality at the quantum level}

In this subsection we prove the theorem from representation theory that is
predicted by
the equality (of the numbers of lattice points) of the momentum polyhedra for
bending
and Gel'fand-Tsetlin. The required facts from representation theory can be
found in \cite{FultonHarris} and \cite{Zhelobenko}.

\begin{theorem}\label{quantumduality}
The multiplicity of the $1$-dimensional representation $(det)^{\vL}$ in
$\otimes_1^n \mathcal{S}^{r_j}(V) $ is equal to the multiplicity of
the weight $\br$ in the
irreducible representation $C^{\vL}\bigwedge^{m+1}V$ of $\Un$. This common
multiplicity is in fact equal to the number of lattice points in ${\bf P}$.
\end{theorem}

The theorem will be a consequence of the next three lemmas. We will need

\begin{definition} Let $\lambda$ be an $l$-tuple of positive integers
and $\mu$ be a partition. Then the Kostka number
$K_{\mu \lambda}$ is the number of ways to fill in the Young diagram
corresponding
to $\mu$ with $\lambda_1$ $1$'s, $\lambda_2$ $2$'s, $\ldots$, $\lambda_l$ $l$'s
so
that the rows are weakly increasing and the columns are strongly increasing.
\end{definition}

By applying Pieri's formula iteratively one gets  \cite[(A.9)]{FultonHarris}:
\begin{lemma}
$$\mathcal{S}^{r_1}(V) \otimes \mathcal{S}^{r_2}(V) \otimes \cdots \otimes
 \mathcal{S}^{r_n}(V)=\oplus_{\mu} K_{\mu \br}V(\mu).$$
\end{lemma}

We obtain
\begin{corollary}
The multiplicity of the $1$-dimensional representation $(det)^{\vL}$ in
$\otimes_1^n \mathcal{S}^{r_j}(V) $ is equal to the Kostka number $K_{\vL(1^{m+1})\br}$.
\end{corollary}

Here the symbol $\vL(1^{m+1})$ means the partition
$(\vL,\vL,\ldots,\vL)$ (there are $m+1$ $\vL$'s).
The corresponding Young diagram has $m+1$ rows and $\vL$ columns.

In order to compare $K_{\vL(1^{m+1})\br}$ with the multiplicity of the weight
$\br$ in the irreducible representation $C^{\vL}\bigwedge^{m+1}V$ of $\Un$ we
recall
there is a basis for an irreducible representation of $GL(n)$ labeled by
semistandard
Young tableaux. Suppose the highest weight of the representation is $\mu$. We
also use
$\mu$ to denote the Young diagram associated to $\mu$. A {\em semistandard
filling} of
$\mu$ is an assignment of the integers between $1$ and $n$ to the boxes of
$\mu$
such that the rows are weakly increasing and the columns are strongly
increasing.
The associated basis is a weight basis, and the weight of of the basis vector
corresponding to a semistandard tableau is $(k_1,\ldots,k_n)$, where $k_i$ is
the
number of $i$'s in the tableau. Thus we have proved

\begin{lemma}
$K_{\vL(1^{m+1})\br}$ is also the multiplicity of the weight $\br$ in
$C^{\vL}\wedge^{m+1}V$ of $GL(n)$.
\end{lemma}

It still remains to prove that the number of lattice points in ${\bf P}$ is the
common
multiplicity.

To see this we recall that there is an orthonormal basis (the Gel'fand-Tsetlin
basis)
for
the irreducible representation $C^{\vL}\wedge^{m+1}V$ indexed by
Gel'fand-Tsetlin
patterns whose top row consists of $m+1$ $\vL$'s and
$n-m-1$ zeroes. Moreover, this basis is a weight basis, and the weight
of a basis vector corresponding to a Gel'fand-Tsetlin pattern is given by the
differences
in the row sums starting with the bottom entry in the pattern. Thus we have

\begin{lemma}
The number of lattice points in ${\bf P}$ is equal to the dimension of the
$\br^{\text{th}}$ weight space in $C^{\vL}\bigwedge^{m+1}V$.
\end{lemma}

It follows that the count of lattice points in ${\bf P}$ gives the correct
answer for both multiplicities.

\begin{remark}
One might ask whether there is a direct combinatorial argument to establish
the last lemma above, i.e. that the number of semistandard Young tableaux of
weight $\br$
is equal to the number of Gel'fand-Tsetlin patterns of weight $\br$. In fact,
there is a
one to one weight preserving correspondence between semistandard
Young tableaux and Gel'fand-Tsetlin patterns, see \cite{Gel'fandZelevinsky}.
\end{remark}

\section{Appendix: Bending and Hitchin Hamiltonians}
In \cite{AdamsHarnadHurtubise}, the authors constructed a  duality
between integrable systems on certain pairs of finite dimensional coadjoint
orbits of loop groups. The spaces considered in Chapter 8 above belong to their
family.  However we shall show below that the Hamiltonians considered
in \cite{AdamsHarnadHurtubise} are different from ours even for the
case of polygons in $\mathbb{R}^3$ i.e. $m=1$. The point is that
our Hamiltonians depend on a triangulation of the polygon by diagonals
(equivalently, the Gelfand-Tsetlin decomposition into increasing
principal diagonal blocks). This introduces an asymmetry which is not
present in the theory of \cite{AdamsHarnadHurtubise}.
We will not review the theory of \cite{AdamsHarnadHurtubise} here but
will instead review how one obtains integrable systems on $\Mr$
by associating a matrix $A(z)$ with entries which are polynomials in a
complex variable $z$ to a point $\be \in \Mrtw$. The construction of
\cite{AdamsHarnadHurtubise} is more general but reduces to the above
for the case in hand.

We will
refer to the resulting Hamiltonians as Hitchin Hamiltonians although
the construction we are about to describe antedated and is a very special
case of Hitchin's construction
of integrable systems associated to Higgs fields.  We will follow
the notation of \cite{Hitchin}, pages 46--52.

\subsection{Hitchin Hamiltonians}
We will describe the construction of Hitchin Hamiltonians for the
case considered in this paper. So let $\be \in \Mrtw$. Choose $n$ points
$\alpha_1,\alpha_2,\cdots,\alpha_n$. Put $p(z) =
\prod\limits_{i=1}^n (z-\alpha_i)$ and define a
matrix--valued polynomial $A(z)$ by
$$A(z) = p(z) \sum\limits_{i=1}^n \frac{e_i}{z -\alpha_i}.$$
We then consider the characteristic polynomial
$\det(w\id - A(z))$ of $A(z)$. The coefficients of the characteristic
polynomial (considered as functions on $\Mr$) are the Hitchin
Hamiltonians. In the case in which $m=1$ it suffices to consider the
functions
$$H_j(\be) = \sum\limits_{i \neq j} \frac{e_i\cdot e_j}{\alpha_i - \alpha_j},
1 \leq j \leq n.$$
We will now examine when the bending flows leave invariant the Hitchin
Hamiltonians for the special case $m=1$. We note that if $f$ is a smooth
function of $\Mr$ in order to compute its derivative along a bending flow
we may choose any lift of $f$ to $\Mrtw$ and any lift of the bending
flow to $\Mrtw$.

\subsection{Bending Hamiltonians are not Hitchin Hamiltonians}

In this section we will prove
\begin{proposition}
If $\alpha_1,\cdots,\alpha_n$ are distinct and $n \geq 5$, then the diagonal
$d_{13}:= ||e_1+e_2||$ does not belong to the Hitchin system.
\end{proposition}

The proposition will follow from the next two lemmas. Let $\phi_t$
be the bending flow associated to diagonal $d_{13}$. We will lift $d_{13}$
to $\Mrtw$ and lift $\phi_t$
to $\Mrtw$ so that the edges $e_1$ and $e_2$ are rotated around the
diagonal $e_1 + e_2$ and the remaining edges are left fixed. We will
continue to use $\phi_t$ for the lifted flow.

\begin{lemma}
Let $i > 2$ and $c$ be a constant. Then the function $g(\be) :\  = c \ e_1\cdot
e_i$ is invariant under
the flow $\phi_t
\Leftrightarrow c=0$.
\end{lemma}
\begin{proof}
We have
$$\dot{e_1} = e_2 \times e_1, \  \dot{e_2} = e_1 \times e_2 \ \text{and} \
\dot{e_i} = 0, \  i > 2.$$
Hence, if $i>2$ we have
$$\frac{d}{dt}\, c\, e_1\cdot e_i = c\,(e_2 \times e_1)\cdot e_i \equiv 0.$$
But there exist polygons such that the scalar triple $(e_2 \times e_1)\cdot
e_i$
is nonzero. Consequently $c=0$.
\end{proof}

In the proof of the next two lemmas we will use the following notation.
If $f$ and $g$ are two smooth functions on $\Mrtw$ we will write
$f \equiv g$ if $f$ and $g$ differ by a function that is invariant under
$\phi_t$.

Recall $H_n(e)$ is the Hitchin Hamiltonian defined by
$$H_n(e) = \sum\limits_{i=1}^{n-1} \frac{e_i\cdot e_n}{\alpha_i - \alpha_n}.$$

\begin{lemma}\label{littleflap}
$H_n$ is invariant under $\phi_t \Leftrightarrow \alpha_1 = \alpha_2$.
\end{lemma}

\begin{proof}

We have, since $e_i\cdot e_n$ is invariant under $\phi_t$ for $i > 2$
(since $e_i$ is)
\begin{align*}
H_n(e) \equiv & \frac{e_1\cdot e_n} {\alpha_1 - \alpha_n} +
\frac{e_2\cdot e_n} {\alpha_2 - \alpha_n}\\
= & \frac{e_1\cdot e_n} {\alpha_1 - \alpha_n} +
\frac{((e_1+e_2) - e_1) \cdot e_n} {\alpha_2 - \alpha_n}\\
\equiv & \ ( \frac{1}{\alpha_1 -\alpha_n} - \frac{1}{\alpha_2 -\alpha_n}) \
e_1\cdot e_n.
\end{align*}
The last line follows because $e_1 + e_2$ is invariant under $\phi_t$.

By the previous lemma we find that $H_n(e)$ is invariant if and only if
$$\frac{1}{\alpha_1 -\alpha_n} - \frac{1}{\alpha_2 -\alpha_n}=0$$
\end{proof}

There remains the possibility that all the bending Hamiltonians could
be Hitchin Hamiltonians for a singular curve. We will now show that this
cannot happen (for the case $n=5$). Suppose then $n=5$. Let $\phi_t$ be the
bending flow
along the diagonal $e_1 + e_2$ and $\psi_t$ be the bending flow
along the diagonal $e_1 + e_2 + e_3$. We will lift $\phi_t$ (resp. $\psi_t$)
to $\Mrtw$ so that only the first two (resp. three) edges are moved. We will
now use the notation
$f\equiv g$ to denote that $f$ and $g$ differ by a function that
is invariant under both flows.

\begin{lemma}\label{bigflap}
Suppose the Hitchin Hamiltonian $H_5$ is invariant under
both bending flows $\phi_t$ and $\psi_t$.
Then the first three $\alpha_i$'s are equal.
\end{lemma}

\begin{proof}
Since $e_4$ and $e_5$ are invariant under both flows we have
\begin{align*}
H_5(e) \equiv & \frac{e_1\cdot e_5} {\alpha_1 - \alpha_5} +
\frac{e_2\cdot e_5} {\alpha_2 - \alpha_5}
+ \frac{e_3\cdot e_5} {\alpha_3 - \alpha_5}\\
= & \frac{e_1\cdot e_5} {\alpha_1 - \alpha_5} +
\frac{e_2\cdot e_5} {\alpha_2 - \alpha_5} +
\frac{((e_1+e_2 + e_3) - (e_1+e_2)) \cdot e_5} {\alpha_3 - \alpha_5}\\
\equiv & \ ( \frac{1}{\alpha_1 -\alpha_5} - \frac{1}{\alpha_3 -\alpha_5}) \
e_1\cdot e_5 + ( \frac{1}{\alpha_2 -\alpha_5} - \frac{1}{\alpha_3 -\alpha_5})
e_2\cdot e_5.
\end{align*}

But by Lemma \ref{littleflap} we have
$$\frac{1}{\alpha_1 -\alpha_5} = \frac{1}{\alpha_2 -\alpha_5}.$$

Consequently we have
$$H_5(e) \equiv ( \frac{1}{\alpha_1 -\alpha_5} - \frac{1}{\alpha_3 -\alpha_5})
(e_1 + e_2) \cdot e_5.$$
We put $c$ equal to the coefficient of $(e_1 + e_2) \cdot e_5$ act by $\psi_t$
and differentiate to obtain
$$\frac{d}{dt} c (e_1 + e_2)\cdot e_5 =
c((e_1 + e_2 + e_3)\times (e_1 + e_2 ))\cdot e_5 = 0.$$
But there exist polygons such that the number
$((e_1 + e_2 + e_3)\times (e_1 + e_2 ))\cdot e_5$
is nonzero. Consequently $c=0$.

Hence $\alpha_1 = \alpha_3$ and the first three $\alpha_i$'s are equal.

\end{proof}

Now we prove that if $H_1(e)$ is also invariant under $\phi_t$ and
$\psi_t$ then all the $\alpha_i$'s are equal.

\begin{lemma}
Suppose the Hitchin Hamiltonian $H_1$ is invariant under
both bending flows $\phi_t$ and $\psi_t$.
Then the last three $\alpha_i$'s are equal.
\end{lemma}

\begin{proof}
Lift $\psi_t$ from $\Mr$ to $\Mrtw$ so that it rotates the triangle
with edges $e_1 + e_2 + e_3, e_4$ and $e_5$ around the diagonal
$e_1 + e_2 + e_3$ (so $\psi_t$ leaves $e_1,e_2$ and $e_3$ fixed). Repeat the
argument
in Lemma \ref{littleflap}
to find $\alpha_4 = \alpha_5$. Now lift $\phi_t$ to $\Mrtw$ so
that it rotates the quadrilateral with edges $e_5,e_4,e_3$ and $e_1 + e_2$
around the diagonal $e_1 + e_2$ (so $\phi_t$ leaves $e_1$ and $e_2$ fixed).
Repeat the argument of Lemma \ref{bigflap}
to find the last three $\alpha_i$'s are equal.
\end{proof}

We can now show for pentagons the bending systems and the Hitchin systems
never coincide.

\begin{theorem}
For the case $n=5$, the Hitchin system never coincides with the bending system.
\end{theorem}
\begin{proof}
If the systems coincide then $H_1$ and $H_5$ are invariant under both bending
flows and
all the $\alpha_i$'s are equal. But in this case we obtain
$$A(z) = \sum\limits_{i=1}^5 \frac{e_i}{z-\alpha_i} = \frac{1}{z-\alpha_1}
\sum\limits_{i=1}^5 e_i = 0.$$
Hence $\det(w \id - A(z))$ = $w^2$ and there are no nontrivial Hitchin
Hamiltonians.
\end{proof}

\end{document}